\documentclass[aop,MSNbibl,seceqn,citesort,dvips]{arximspdf}

\doi{10.1214/11-AOP692}
\volume{41}
\issue{3B}
\pubyear{2013}
\firstpage{1978}
\lastpage{2012}

\makeatletter

\renewcommand{\diagup}{/}

\newtheorem{theorem}{Theorem}[section]
\newproclaim{condition}{Hypothesis}

\newproclaim{definition}{Definition}[section]
\newproclaim{example}{Example}
\newtheorem{lemma}{Lemma}[section]
\newproclaim{remark}{Remark}[section]
\newproclaim{notation}{Notation}[section]

\makeatother

\begin{document}
\begin{frontmatter}

\title{Regularity of solutions to quantum master equations: A stochastic approach}
\runtitle{Regularity of solutions to QMEs}

\begin{aug}
\author[A]{\fnms{Carlos M.} \snm{Mora}\corref{}\thanksref{t1}\ead[label=e1]{cmora@ing-mat.udec.cl}\ead[label=u1,url]{www.ing-mat.udec.cl/\textasciitilde cmora/}}
\runauthor{C. M. Mora}
\affiliation{Universidad de Concepci\'{o}n}
\address[A]{CI$^2$MA,
Departamento de Ingenier\'{\i}a Matem\'{a}tica \\
Facultad de Ciencias F\'{\i}sicas y Matem\'{a}ticas\\
Universidad de Concepci\'{o}n \\
Casilla 160 C, Concepci\'{o}n\\
Chile\\
\printead{e1}\\
\printead{u1}} 
\end{aug}

\thankstext{t1}{Supported in part by FONDECYT Grants 1070686 and
1110787, and by BASAL Grants PFB-03 and FBO-16 as well as by PBCT-ACT
13 project.}

\received{\smonth{3} \syear{2011}}
\revised{\smonth{7} \syear{2011}}

%
\begin{abstract}
Applying probabilistic techniques we study regularity properties of
quantum master equations (QMEs) in the Lindblad form with unbounded
coefficients; a density operator is regular if, roughly speaking, it
describes a quantum state with finite energy. Using the linear
stochastic Schr\"{o}dinger equation we deduce that solutions of QMEs
preserve the regularity of the initial states under a general
nonexplosion condition. To this end, we develop the probabilistic
representation of QMEs, and we prove the uniqueness of solutions for
adjoint quantum master equations. By means of the nonlinear stochastic
Schr\"{o}dinger equation, we obtain the existence of regular stationary
solutions for QMEs, under a Lyapunov-type condition.
\end{abstract}

%

\pdfkeywords{60H15, 60H30, 81C20,
46L55, Quantum master equations,
stochastic Schrodinger equations,
regular solutions, probabilistic representations,
open quantum systems}

\end{frontmatter}

\section{Introduction}

In order to establish
the well-posedness of the mean values of quantum observables
represented by unbounded operators,
we investigate
the regularity of solutions of quantum master equations (with unbounded
coefficients) in stationary and transient regimes.
For this purpose,
we use classical stochastic analysis.

\subsection{Gorini--Kossakowski--Lindblad--Sudarshan equations}

In many open quantum systems,
the states of a small quantum system with Hamiltonian $H\dvtx \mathfrak{h}
\rightarrow\mathfrak{h}$
evolve according to the operator equation
%
%
\begin{equation}
\label{3}
\frac{d}{dt} \rho_{t}( \varrho) = \mathcal{L}_{\ast}(
\rho_{t}( \varrho) ),\qquad \rho_{0}( \varrho) =\varrho,
\end{equation}
where
$
\mathcal{L}_{\ast}( \rho)
=
G \rho+ \rho G^{\ast} +\sum_{k=1}^{\infty}L_{k} \rho L_{k}^{\ast}
$
(see, e.g.,~\cite{BreuerPetruccione2002,GardinerZoller2004,WisemanMilburn2010}).
Here,
$( \mathfrak{h}, \langle\cdot,\cdot\rangle)$ is a separable complex
Hilbert space,
$G, L_{1}, L_{2}, \ldots$ are given linear operators in $\mathfrak{h}$
satisfying
$
G=-iH-\frac{1}{2} \sum_{k=1}^{\infty}L_{k}^{\ast}L_{k}
$
on suitable common domain
and
the unknown density operator $\rho_{t}( \varrho)$
is a nonnegative operator in $\mathfrak{h}$ with unit trace.
The operators $L_{1},L_{2},\ldots$ describe the weak interaction between
the small quantum system and a heat bath.

The measurable physical quantities of the small quantum system are
represented by self-adjoint operators in $\mathfrak{h}$, which are
called observables.
Very important observables are unbounded,
like position, momentum and kinetic energy operators.
In the Schr\"odinger picture,
the mean value of the observable $A$ at time $t$ is given by
$\operatorname{tr} ( \rho_{t}(\varrho) A )$,
the trace of $\rho_{t}( \varrho) A$.

In the Heisenberg picture,
the initial density operator $\varrho$ is fixed. Using,
for instance, (\ref{3})
we obtain the following equation of motion for the observable~$A$:
%
%
\begin{equation}
\label{41}
\frac{d}{dt} \mathcal{T}_{t}( A )
=
\mathcal{L} ( \mathcal{T}_{t}( A ) ),\qquad
\mathcal{T}_{0}( A ) = A,
\end{equation}
where
$
\mathcal{L} ( \mathcal{T}_{t}( A ) )
=
\mathcal{T}_{t}( A ) G
+ G^{\ast} \mathcal{T}_{t}( A ) +\sum_{k=1}^{\infty} L_{k}^{\ast}
\mathcal{T}_{t}( A ) L_{k}
$; see, for example, \mbox{\cite{BreuerPetruccione2002,GardinerZoller2004}}.
The expected value of $A$ at time time $t$
is given by $\operatorname{tr} (\varrho\mathcal{T}_{t}( A ) )$.

\subsection{\texorpdfstring{Stochastic Schr\"{o}dinger equations
(SSEs)}{Stochastic Schrodinger equations (SSEs)}}

The evolution of the state of a quantum system conditioned on
continuous measurement
is governed (see, e.g., \cite
{BarchielliGregoratti2009,Belavkin1989,WisemanMilburn2010}) by
the stochastic evolution equation on $\mathfrak{h}$.
%
%
\begin{equation}
\label{5}
Y_{t}=Y_{0}+\int_{0}^{t}G( Y_{s}) \,ds+\sum_{k=1}^{\infty}\int
_{0}^{t}L_{k}( Y_{s}) \,dB_{s}^{k}.
\end{equation}
Here $ G( y) = G y + \sum_{k=1}^{\infty}(\Re\langle y,L_{k} y \rangle
L_{k} y -\frac{1}{2} \Re^{2}\langle y,L_{k} y \rangle y) $, $ L_{k}( y)
=L_{k} y -\break \Re\langle y, L_{k} y \rangle y $ and $ B^{1}, B^{2}, \ldots$
are real valued independent Wiener processes.

\begin{example}
\label{exmeasurement}
Set $\mathfrak{h} = L^{2}( \mathbb{R},\mathbb{C})$.
Let $Q,P\dvtx \mathfrak{h} \rightarrow\mathfrak{h}$ be defined by
$ Q f ( x ) = x f( x ) $
and
$ P f ( x ) = -i f ^{\prime} ( x ) $.
In (\ref{5}),
take
$H = \frac{1}{2m} P^2 + c Q^2$,
$ L_{1} = \alpha Q $
and
$ L_{2} = \beta P $,
with
$m>0$,
$\alpha, \beta\geq0$
and $c \in\mathbb{R}$.
For all $k\geq3$,
fix
$L_{k}=0$.
\end{example}

Example~\ref{exmeasurement} with $\alpha, \beta, c >0$
describes the simultaneous monitoring of position and momentum of
a linear harmonic oscillator; see, for example, \cite
{GoughSobolev2004,ScottMilburn2001}.
Taking instead
$\alpha>0$ and $\beta= c = 0$
we get a well-studied model
for the continuous measurement of position
of a free particle; see, for example, \cite
{Diosi1988,GoughSobolev2004,BassiDurrKolb2010,Kolokoltsov2000}
and references therein.

Our main tool for studying (\ref{3}) and (\ref{41}) is the following linear SSE
on~$\mathfrak{h}$:
%
%
\begin{equation}
\label{2}
X_{t}( \xi) =
\xi+\int_{0}^{t}GX_{s}( \xi) \,ds + \sum_{k=1}^{\infty}\int
_{0}^{t}L_{k}X_{s}( \xi) \,dW_{s}^{k},
\end{equation}
where
$ W^{1}, W^{2}, \ldots$ are real valued independent Wiener processes
on a filtered complete probability space
$( \Omega,\mathfrak{F},(\mathfrak{F}_{t}) _{t\geq0},\mathbb{P}) $.
In fact,
the basic assumption of this paper is:
\begin{longlist}[(H)]
\item[(H)]
There exists a nonnegative self-adjoint operator $C$ in $\mathfrak{h}$
such that:
(i)
$G$ is relatively bounded with respect to $C$;
and
(ii)
(\ref{2}) has a unique $C$-solution
for any initial condition $\xi$ satisfying
$ \mathbb{E} ( \Vert C \xi\Vert^{2} + \Vert\xi\Vert^{2} ) <
\infty$.
\end{longlist}
Here,
a strong solution $X ( \xi)$ of (\ref{2}) is called $C$-solution
if
$\mathbb{E}\Vert X_{t}( \xi) \Vert^{2}\leq\mathbb{E}\Vert\xi\Vert^{2}$,
and
the function
$t \mapsto\mathbb{E}\Vert C X_{t}( \xi) \Vert^{2}$
is uniformly bounded on compact time intervals; see Definition \ref
{definicion2} for details.

The law of $X_{t}( Y_0 ) / \| X_{t}( Y_0 ) \|$
with respect to
$\Vert X_{T}( Y_0 ) \Vert^{2}\cdot\mathbb{P}$
coincides with the law of $Y_{t}$
for all $t \in[0, T]$; see~\cite{MoraReAAP2008}.
The main technical and conceptual advantage of (\ref{5}) over (\ref{2})
is that the norm of $Y_t$ is equal to $1$.

\subsection{Principal objectives}
\label{subsecobjectives}

Our main goal is to make progress in the understanding of the evolution of
$\operatorname{tr} ( \rho_{t}(\varrho) A )$
when $A$ is unbounded.

Given a self-adjoint nonnegative operator $C$ in $\mathfrak{h}$,
we denote by $\mathfrak{L}_{1,C}^{+}( \mathfrak{h}) $
the set of all nonnegative operators $\varrho\dvtx \mathfrak{h}
\rightarrow\mathfrak{h}$ for which, loosely speaking,
$C \varrho$ is a trace-class operator; see Definition~\ref{def2}.
From Section~\ref{secprob-rep} we have that
the expected value of $A$ with respect to
$\varrho\in\mathfrak{L}_{1,C}^{+}( \mathfrak{h}) $ is well defined
whenever
$A \in\mathfrak{L}( ( \mathcal{D}( C)$, $\mbox{$\Vert\cdot\Vert_{C}$})
,\mathfrak{h})$,
where $\mathfrak{L}( ( \mathcal{D}( C), \mbox{$\Vert\cdot\Vert_{C}$})
,\mathfrak{h})$
is the space of all operators relatively bounded with respect to $C$.
Our first objective is:
\begin{longlist}[(O1)]
\item[(O1)] To prove that the solution
$\rho_{t}( \varrho)$ of (\ref{3}) belongs to $\mathfrak{L}_{1,C}^{+}(
\mathfrak{h})$
(for all $t > 0$)
provided that
$C$ satisfies hypothesis (H)
and that $ \varrho\in\mathfrak{L}_{1,C}^{+}( \mathfrak{h})$.
\end{longlist}

The key condition to guarantee
the uniqueness of solution of (\ref{3})
is the existence of a self-adjoint nonnegative operator $C$ in
$\mathfrak{h}$ such that formally
%
%
\begin{equation}
\label{i2}
\mathcal{L} ( C^2 ) \leq K ( C^2 + I ),
\end{equation}
where $I$ is the identity operator in $\mathfrak{h}$ and $K \in[ 0,
\infty[$.
This condition,
introduced by Chebotarev and Fagnola~\cite{ChebFagn1998}
(see also~\cite{Chebotarev2000,Fagnola1999,GarQue1998}),
is a quantum analog of the Lyapunov condition for nonexplosion of
classical Markov processes; see~\cite{Chebotarev2000} for heuristic arguments.
Since
hypothesis (H) holds under a weak version of (\ref{i2})
(see~\cite{FagnolaMora2010,MoraReIDAQP2007} and Remark \ref
{notaSuffCond}), inequality (\ref{i2}) is the underlying assumption of
objective (O1).
In many physical examples,
relevant observables belong to
$\mathfrak{L}( ( \mathcal{D}( C), \mbox{$\Vert\cdot\Vert_{C}$}) ,\mathfrak{h})$
for\vspace*{1pt} some $C$ satisfying (\ref{i2}).
In Example~\ref{exmeasurement},
for instance,
$C = P^2 + Q^2$ satisfies hypothesis (H)
(see, e.g.,~\cite{FagnolaMora2010,MoraReIDAQP2007}),
and
the position and momentum operators
$Q$ and $P$
are $( P^2 + Q^2 )$-bounded.

Previously,
the regularity of the solutions to (\ref{3})
has been treated in~\cite{ArnoldSparber2004,ChebGarQue1998,Davies1977b}
using methods from the operator theory.
Exploiting the characteristics of a model describing
a variable number of neutrons moving in a translation invariant
external reservoir of unstable atoms,
Davies~\cite{Davies1977b} established that
$\rho_{t}(\varrho) \in\mathfrak{L}_{1,C}^{+}( \mathfrak{h})$
whenever $C$ is the particle number operator on an adequate Fermion
Fock space.
Arnold and Sparber~\cite{ArnoldSparber2004} obtained the same property with
$C$ being essentially the energy operator
for a linear quantum master equation associated to a diffusion model
with Hartree interaction.\vadjust{\goodbreak}

The second objective presents
the first attempt (to the best of my knowledge)
to show the existence of stationary solutions of (\ref{3}) with finite energy.
\begin{itemize}[(O2)]
\item[(O2)] To prove the existence of
a stationary solution of (\ref{3}) belonging to $\mathfrak
{L}_{1,D}^{+}( \mathfrak{h})$
provided essentially that:
\begin{itemize}[(L)]
\label{PageHypL}
\item[(L)] There exist two nonnegative self-adjoint operators $C$ and $D$
and a constant $K>0$
such that
$\{ x \in\mathfrak{h}\dvtx \| D x \| ^2 + \| x \| ^2 \leq1 \}$
is compact in $\mathfrak{h}$,
and
$
\mathcal{L} ( C^2 ) \leq-D^2 + K ( 1 + I )
$.
\end{itemize}
\end{itemize}

Hypothesis (L) is a quantum version of
the Lyapunov criterion for the existence of invariant probability
measures for stochastic differential equations,
which applies to many open quantum systems; see, for example, \cite
{FagReb2001} and Section~\ref{secoscillator}.
Let
$\mathfrak{L}( \mathfrak{h} ) $
be the set of all bounded operators from $\mathfrak{h}$ to $\mathfrak{h}$.
In the case where
$G$ is the infinitesimal generator of a strongly continuous contraction
semigroup on $\mathfrak{h}$
and, loosely speaking,
$ \mathcal{T}_{t} ( I ) = I$ for all $t \geq0$,
Fagnola and Rebolledo~\cite{FagReb2001} proved that under hypothesis (L),
there exists at least one density operator
$\varrho_{\infty}$
satisfying
$
\operatorname{tr} (\varrho_{\infty} \mathcal{T}_{t}( A ) )
=
\operatorname{tr} (\varrho_{\infty} A )
$
for all $t \geq0$ and $A \in\mathfrak{L} ( \mathfrak{h})$.
The main point of objective (O2)
is that
among such stationary states $\varrho_{\infty} $
we can select a finite-energy density operator belonging to the domain
of $ \mathcal{L}_{\ast}$,
under the same hypothesis (L).

The third main objective develops
the rigorous probabilistic representation of solution of (\ref{3}),
the key step to achieve objectives (O1) and (O2).
\begin{longlist}[(O3)]
\item[(O3)]
Assume hypothesis (H),
and let
$
\varrho= \mathbb{E} \vert\xi\rangle\langle\xi\vert
$,
where $\xi$ is a $\mathfrak{h}$-valued random variable such that
$
\mathbb{E} \Vert\xi\Vert^{2} = 1
$
and
$
\mathbb{E} \Vert C \xi\Vert^{2} < \infty
$.
We wish to prove that (\ref{3}) has a unique solution,
which is
%
%
\begin{equation}
\label{I1}
\rho_{t}( \varrho)
= \mathbb{E} \vert X_{t}( \xi) \rangle\langle X_{t}( \xi) \vert.
\end{equation}
\end{longlist}
In Dirac notation,
$\vert x\rangle\langle y\vert\dvtx \mathfrak{h} \rightarrow\mathfrak{h}$
is defined by
$
\vert x\rangle\langle y\vert( z ) = \langle y,z\rangle x
$,
with \mbox{$x,y \in\mathfrak{h}$}.
Using (\ref{I1}) we can assert that
%
%
\begin{equation}
\label{i9}
\rho_{t}( \varrho) =
\mathbb{E} \vert Y_{t} \rangle\langle Y_{t} \vert
\end{equation}
with $Y_0 = \xi$
(see~\cite{MoraReAAP2008}).
Objective (O3), together with (\ref{i9}),
shows that
physical models based on the stochastic Schr\"odinger equations
are in good agreement with
their formulations in terms of
quantum master equations.

In the physical literature,
the probabilistic representations
(\ref{I1}) and (\ref{i9})
of the density operator at time $t$
have been obtained by means of formal computations; see, for example,
\cite{BarchielliBelavkin1991,BreuerPetruccione2002,GisinPercival1992}.
Barchielli and Holevo~\cite{BarchielliHolevo1995} established
essentially (\ref{I1}) and (\ref{i9}) in situations where
$G,L_{1},L_{2},\ldots$ are bounded.

\subsection{Approach}
\label{subsecapproach}

In the perspective of the operator theory,
methods based on the Hille--Yosida theorem
and perturbations of linear operators~\cite{Kato1995,Pazy1983}
present severe limitations for studying
linear functionals of the solutions of (\ref{3}) and (\ref{41}).
For example,
it is very difficult
to decompose
$ \mathcal{L}_{\ast}$
into
$ \mathcal{L}_{\ast}^1+ \mathcal{L}_{\ast}^2$
for
a dissipative operator $ \mathcal{L}_{\ast}^1$ in $\mathfrak{L}_{1}(
\mathfrak{h})$
and
an infinitesimal generator $ \mathcal{L}_{\ast}^2$ of a $C_0$ semigroup
of contractions on\vadjust{\goodbreak}
$\mathfrak{L}_{1}( \mathfrak{h})$,
which together satisfy
$
\Vert\mathcal{L}_{\ast}^1 ( \varrho) \Vert_{1}
\leq
\alpha\Vert\mathcal{L}_{\ast}^2 ( \varrho) \Vert_{1} + K \Vert
\varrho\Vert_{1}
$
whenever $\varrho\in\mathcal{D} ( \mathcal{L}_{\ast}^2 )$.
Here,
$0 \leq\alpha< 1$,
$K\geq0$
and
$\mathfrak{L}_{1}( \mathfrak{h})$ is
the Banach space of trace-class operators on $\mathfrak{h}$ equipped
with the trace norm
\mbox{$ \Vert\cdot\Vert_{1}$}.
Another difficulty is that
$ \mathcal{L}$ and $ \mathcal{L}_{\ast}$
are defined formally;
indeed
$ \mathcal{L}$ and $ \mathcal{L}_{\ast}$
can be interpreted as sesquilinear forms,
but without having a priori knowledge about their cores.

When $G$ is the generator of a $C_0$ semigroup of contractions on
$\mathfrak{h}$, Davis~\cite{Davies1977a} provided solutions of (\ref
{3}) by means of semigroups.
Modifying Davis's ideas,
Chebotarev constructed a quantum dynamical semigroup $\mathcal{T}^{(\min)}$
that is weak solution of (\ref{41})
by generalizing the Chung construction of the minimal solution of
Feller--Kolmogorov equations for countable state Markov chains; see
Remark~\ref{nota7}.
Under certain conditions involving (\ref{i2}) and invariant sets for
$\exp(Gt )$,
Chebotarev and Fagnola~\cite{ChebFagn1998} proved the uniqueness of
$\mathcal{T}^{(\min)}$; see Remark~\ref{nota7}.
This property implies that
$ \mathcal{L}_{\ast}$ is the infinitesimal generator of the predual
semigroup $\rho^{(\min)}$ of $\mathcal{T}^{(\min)}$,
and a core for $ \mathcal{L}_{\ast}$ is formed by the linear span of
all $\vert x\rangle\langle y\vert$
with $x,y$ belonging to $\mathcal{D} ( G )$,
the domain of $G$;
see Remark~\ref{nota8}.
In Remark~\ref{nota8},
we outline how to obtain
$\rho^{(\min)} ( \mathfrak{L}_{1,C}^{+}( \mathfrak{h}) ) \subset
\mathfrak{L}_{1,C}^{+}( \mathfrak{h}) $
under various assumptions including
$ \exp( G t ) \mathcal{D}( C ) \subset\mathcal{D}( C )$.
It is a hard problem, in general, to find $C$ satisfying (\ref{i2})
whose domain $ \mathcal{D}( C )$ is invariant under the action of $
\exp( G t ) $.

In contrast to closed quantum systems,
solutions of (\ref{3}) are not decomposable as dyadic products of
solutions of evolution equations in $\mathfrak{h}$.
Nevertheless,
the solution of (\ref{3}) is unraveled into stochastic quantum trajectories;
more precisely,
objective (O3) establishes that
$
\rho_{t}( \varrho)
$
is expressed as the mean value of quadratic functionals of the
solutions of SSEs
in a general context.
This property allows us to achieve objectives (O1) and (O2) by using SSEs,
without serious difficulties
and
without assumptions involving invariant sets for $\exp(Gt )$.
Applying (\ref{I1})
we also deduce that
$
\rho_{t}( \varrho)
$
satisfies (\ref{3}) in both sense integral and $\mathfrak{L}_{1}(
\mathfrak{h}) $-weak.
This leads to prove rigorously
some dynamical properties of $ \rho_{t}( \varrho)$ given in physics;
see, for example, Theorem~\ref{teorema5}.

We now focus on objective (O1).
By Section~\ref{secprob-rep},
$ \varrho\in\mathfrak{L}_{1,C}^{+}( \mathfrak{h})$
iff
there exists a $\mathfrak{h}$-valued random variable $\xi$ satisfying
$
\mathbb{E} ( \Vert C \xi\Vert^{2} + \Vert\xi\Vert^{2} ) < \infty
$
and
$
\varrho= \mathbb{E} \vert\xi\rangle\langle\xi\vert
$.
Therefore
(\ref{I1}) leads directly to objective (O1)
since
$
\mathbb{E} \Vert C X_{t}( \xi) \Vert^{2}
+
\mathbb{E} \Vert X_{t}( \xi) \Vert^{2}
< \infty
$.
Assumption (\ref{i2}) is natural
in the context of (\ref{2})
because (\ref{i2}) is essentially the dissipative condition for (\ref{2}).

We turn to objective (O2).
Here,
hypothesis (L) is a classical Lyapunov condition for (\ref{2}).
Relation (\ref{I1}) suggests us that
$ \int_{\mathfrak{h}}\vert x\rangle\langle x\vert\mu( dx)$
is a good candidate for being a stationary solution for (\ref{3})
when
$\mu$ is an invariant probability measure for (\ref{2})
such that $ \int_{\mathfrak{h}} \| x \|^2 \mu( dx) = 1$.
This reduces objective (O2) to prove that
there exists an invariant probability measure for (\ref{2}),
different from the Dirac measure at $0$,
which is a difficult problem.
We instead use (\ref{i9}).
Under a weak version of hypothesis (L),
there exists an invariant probability measure $\Gamma$ for (\ref{5})
such that
$ \int_{\mathfrak{h}} \| x \|^2 \Gamma( dx) = 1$
and
$ \int_{\mathfrak{h}} \| D x \|^2 \Gamma( dx) < \infty$;\vadjust{\goodbreak} see \cite
{MoraReAAP2008}.
Then,
using (\ref{i9}) we deduce that
$
\varrho_{\infty}=\int_{\mathfrak{h}}\vert x\rangle\langle x\vert
\Gamma( dx)
$
is a stationary solution to (\ref{3}) that belongs to
$\mathfrak{L}_{1,D}^{+}( \mathfrak{h})$; see Section~\ref{secStatSol}.
This is a step forward in the study
of the long time behavior of unbounded observables.

\subsection{Technical ideas: Unraveling}
\label{subsecexistence}

Fix $ \varrho\in\mathfrak{L}_{1,C}^{+}( \mathfrak{h})$.
Then
$
\varrho= \mathbb{E} \vert\xi\rangle\langle\xi\vert
$
for some $\mathfrak{h}$-valued random variable $\xi$ satisfying
$
\mathbb{E} ( \Vert C \xi\Vert^{2} + \Vert\xi\Vert^{2} ) < \infty
$; see Section~\ref{secprob-rep}.
We can define
%
%
\begin{equation}
\label{i4}
\rho_{t}( \varrho)
: = \mathbb{E} \vert X_{t}( \xi) \rangle\langle X_{t}( \xi) \vert,
\end{equation}
because
$\rho_{t}( \varrho) $ does not depend on the choice of $\xi$; see
Theorem~\ref{teor7}.
We next outline how to establish that
$
\rho_{t}( \varrho)
$
is a solution of (\ref{3}).

Applying It\^o's formula we obtain
\begin{eqnarray*}
\langle X_{t}( \xi) , x \rangle X_{t}( \xi)
&=&
\langle\xi, x \rangle\xi
+ \int_{0}^{t} \bigl(
\langle X_{s}( \xi) ,x \rangle GX_{s}( \xi)
+
\langle G X_{s}( \xi) , x \rangle X_{s}( \xi)
\bigr) \,ds
\\
&&{}
+ \sum_{k=1}^{\infty} \int_{0}^{t} \langle L_k X_{s}( \xi) , x
\rangle L_k X_{s}( \xi)
\,ds
+ M_t
\end{eqnarray*}
with
$
M_t
=
\sum_{k = 1}^{\infty} \int_{0}^{t} (
\langle X_{s}( \xi) , x \rangle L_k X_{s}( \xi)
+
\langle L_k X_{s}( \xi) , x \rangle X_{s}( \xi)
) \,dW^{k}_{s}
$.
Since $M_t$ is a local martingale,
we use stopping times and the dominated convergence theorem
to deduce that
%
%
\begin{eqnarray}
\label{65}
&&\mathbb{E} \langle X_{t}( \xi) , x \rangle X_{t}( \xi)\nonumber\\
&&\qquad=
\mathbb{E} \langle\xi, x \rangle\xi
+ \int_{0}^{t}
\mathbb{E} \langle X_{s}( \xi) ,x \rangle GX_{s}( \xi)
\,ds
\\
&&\qquad\quad{}
+ \int_{0}^{t}
\mathbb{E} \langle G X_{s}( \xi) , x \rangle X_{s}( \xi)
\,ds
+ \sum_{k=1}^{\infty} \int_{0}^{t} \mathbb{E} \langle L_k X_{s}( \xi
) , x \rangle L_k X_{s}( \xi)
\,ds .\nonumber
\end{eqnarray}
Define the operator
$
\mathcal{L}_{*} ( \xi, s ) \dvtx
\mathfrak{h} \rightarrow\mathfrak{h}
$
to be
\[
\mathbb{E} \vert GX_{s}( \xi) \rangle\langle X_{s}( \xi) \vert
+
\mathbb{E} \vert X_{s}( \xi) \rangle\langle G X_{s}( \xi) \vert
+
\sum_{k=0}^{\infty}
\mathbb{E} \vert L_k X_{s}( \xi) \rangle\langle L_k X_{s}( \xi)
\vert.
\]
We now face the major technical difficulties;
we have to prove that
$\mathcal{L}_{*} ( \xi, s )$
is a trace-class operator such that:
(i)
$
\mathcal{L}_{*} ( \xi,t )
=
\mathcal{L}_{\ast}( \rho_{s}( \varrho) )
$;
(ii) the function
$ s \mapsto
\Vert
\mathcal{L}_{*} ( \xi, s )
\Vert_{1}$
is locally bounded;
and
(iii)
$s \mapsto\mathcal{L}_{*} ( \xi, s )$ is weakly continuous in
$\mathfrak{L}_{1}( \mathfrak{h})$.
Then,
applying (\ref{65}) yields
%
%
\begin{equation}
\label{I5}
\rho_{t}( \varrho)
=
\varrho
+
\int_{0}^{t} \mathcal{L}_{\ast}( \rho_{s}( \varrho) ) \,ds,
\end{equation}
where
we understand the integral of (\ref{I5}) in the sense of the
Bochner integral in $\mathfrak{L}_{1}( \mathfrak{h})$.
Thus,
we can deduce that for any $A \in\mathfrak{L}( \mathfrak{h} )$,
%
%
\begin{equation}
\label{I6}
\frac{d}{dt}\operatorname{tr}( A\rho_{t}( \varrho) )
=
\operatorname{tr}( A \mathcal{L}_{\ast}( \rho_{t}( \varrho) ) ) .\vadjust{\goodbreak}
\end{equation}

\subsection{Technical ideas: Uniqueness}
\label{subsecuniqueness}
Recall that $\rho_{t}( \varrho)$ is defined by (\ref{i4})
for any $ \varrho\in\mathfrak{L}_{1,C}^{+}( \mathfrak{h})$.
In order to establish
the uniqueness of the solution of (\ref{3}) under hypothesis (H),
Theorem~\ref{teor8} extends
$\rho_{t}( \varrho)$
to a strongly continuous semigroup $( \rho_{t})_{t\geq0}$
of bounded operators on $\mathfrak{L}_{1} ( \mathfrak{h})$.
Thus,
$( \rho_{t})_{t\geq0}$ belongs to
the class $\mathcal{S}$
formed by all the locally bounded semigroups
$( \widehat{\rho}_{t} )_{t\geq0}$ on $\mathfrak{L}_{1} ( \mathfrak{h})$
such that for any $ x \in\mathcal{D} (C )$:
(i) $t \mapsto\widehat{\rho}_{t} ( \vert x \rangle\langle x \vert) $
is weakly continuous in $\mathfrak{L}_{1} ( \mathfrak{h}))$;
and
(ii) $\widehat{\rho}_{t} ( \vert x \rangle\langle x \vert)$
satisfies (\ref{I6}) in $t=0$; see Theorem~\ref{teorema9} and Lemma \ref
{lema26} for details.
We next outline the proof that
$( \rho_{t})_{t\geq0}$ is the unique element in~$\mathcal{S}$,
and so
(\ref{3}) has a unique solution (in the semigroup sense).

Let
$( \widehat{\rho}_{t}) _{t\geq0} \in\mathcal{S}$.
Taking in mind that $\mathfrak{L} ( \mathfrak{h})$
is the dual of $\mathfrak{L}_{1} ( \mathfrak{h})$,
we consider the semigroup
$ ( \mathcal{T}_{t} ) _{t\geq0} $ on $\mathfrak{L} ( \mathfrak{h})$
which is the adjoint semigroup of $( \widehat{\rho}_{t}) _{t\geq0}$.
Using techniques from operator theory we obtain
in Lemma~\ref{lema20} that
$ ( \mathcal{T}_{t} ) _{t\geq0} $ is a weak solution of (\ref{41}),
namely,
for all $t \geq0 $, $A \in\mathfrak{L} ( \mathfrak{h})$ and
$x \in\mathcal{D} (C ) $ we have
%
%
\begin{equation}
\label{I7}
\frac{d}{dt} \langle x,\mathcal{T}_{t}( A ) x \rangle
=
\langle x, \mathcal{L} ( \mathcal{T}_{t}( A ) ) x \rangle.
\end{equation}

Now, we wish to prove that
$ ( \mathcal{T}_{t} ) _{t\geq0} $ is the unique weak solution of (\ref{41}),
which is an important problem itself; see, for example, \cite
{ChebFagn1993,ChebFagn1998,ChebGarQue1998,Chebotarev2000,Fagnola1999,MoraJFA2008}.
Suppose for a moment that $\mathfrak{h}$ is finite-dimensional
and
$L_k \ne0$ for only a finite number of $k$.
Applying the It\^o formula to
$
\langle X_{s}( x ) , \mathcal{T}_{t -s}( A ) X_{s}( x ) \rangle
$
we deduce that
\begin{eqnarray*}
&&
\langle X_{t}( x ) , A X_{t}( x ) \rangle\\
&&\quad=
\langle x , \mathcal{T}_{t}( A ) x \rangle
+
M_t
\\
&&\qquad{}
+
\int_{0}^{t} \biggl(
\langle X_{s}( x ) ,
\mathcal{L} ( \mathcal{T}_{t-s}( A ) ) X_{s}( x ) \rangle
-
\biggl\langle X_{s}( x ) ,
\frac{d \mathcal{T}_r ( A )}{dr} \bigg| _{r=t-s} X_{s}( x ) \biggr\rangle
\biggr) \,ds
\end{eqnarray*}
with
\[
M_t
=
\sum_{k=1}^{\infty} \int_{0}^{t}
\bigl(
\langle L_{k} X_{s}( x ) , \mathcal{T}_{t-s}( A ) X_{s}( x ) \rangle
+
\langle X_{s}( x ) , \mathcal{T}_{t-s}( A ) L_{k} X_{s}( x ) \rangle
\bigr)\,
dW_s^{k} .
\]
From (\ref{I7}) we obtain
$
\langle X_{t}( x ) , A X_{t}( x ) \rangle
=
\langle x , \mathcal{T}_{t}( A ) x \rangle
+
M_t
$,
and so the martingale property of $M_t$ leads to
$
\mathbb{E} \langle X_{t}( x ) , A X_{t}( x ) \rangle
=
\langle x , \mathcal{T}_{t}( A ) x \rangle
$, and
hence all the elements in $\mathcal{S}$
are the same semigroups,
which implies
$ \widehat{\rho} = \rho$.

In the general case,
$G$ and $L_k$ are unbounded operators.
Therefore
\[
( s, x )
\mapsto
\frac{d}{ds} \langle x,\mathcal{T}_{t-s}( A ) x \rangle
\bigl( =
\langle x, \mathcal{L} ( \mathcal{T}_{t-s}( A ) ) x \rangle
\bigr)
\]
is not continuous on $[ 0, t ] \times\mathfrak{h}$,
and consequently we cannot apply directly
It\^o's formula to
$
\langle X_{s}( x ) , \mathcal{T}_{t -s}( A ) X_{s}( x ) \rangle
$.
We overcome this difficulty in Section~\ref{subsecteorema3}
by applying It\^o's formula to a regularized version of
$\langle x,\mathcal{T}_{t-s}( A ) x \rangle$;
the resulting stochastic integrals
(similar to those in $M_t$)
are only local martingales,
and so we have to use stopping times.

\subsection{Outline}

Section~\ref{secAdQME} addresses the existence and uniqueness of
solutions for the adjoint quantum master equation, as well as its
probabilistic representation.
Section~\ref{secprob-rep} deals with the probabilistic interpretations
of regular density operators.
In Section~\ref{secQME}
we construct Schr\"{o}dinger evolutions by means of stochastic Schr\"
{o}dinger equations
and study the regularity of solutions to (\ref{3}).
Section~\ref{secStatSol} focusses on the existence of regular
stationary solutions for (\ref{3}).
In Section~\ref{secoscillator}
we apply our results to a quantum oscillator.
Section~\ref{secproofs} is devoted to proofs.

\subsection{Notation}
\label{subsecnot}

Throughout this paper,
the scalar product $\langle\cdot,\cdot\rangle$ is linear in the
second variable and anti-linear in the first one. We write $\mathfrak
{B}( \mathfrak{h}) $ for the Borel $\sigma$-algebra on $\mathfrak{h}$.
Suppose that $A$ is a linear operator in $\mathfrak{h}$.
Then $A^{\ast}$ denotes the adjoint of $A$.
If $A$ has a unique bounded extension to $\mathfrak{h}$, then we
continue to write $A$ for the closure of $A$.

Let $\mathfrak{X}$, $\mathfrak{Z}$ be normed spaces.
We write $\mathfrak{L}( \mathfrak{X},\mathfrak{Z}) $ for the set of all
bounded operators from $\mathfrak{X}$ to $\mathfrak{Z}$ (together with
norm \mbox{$ \Vert\cdot\Vert_{\mathfrak{L}( \mathfrak{X},\mathfrak{Z})}$}).
We abbreviate \mbox{$ \Vert\cdot\Vert_{\mathfrak{L}( \mathfrak{X},\mathfrak
{Z})}$} to~\mbox{$ \Vert\cdot\Vert$}, if no misunderstanding is possible,
and
define $\mathfrak{L}( \mathfrak{X}) = \mathfrak{L}( \mathfrak
{X},\mathfrak{X}) $.
By $\mathfrak{L}_{1}^{+}( \mathfrak{h})$ we mean the subset of all
nonnegative trace-class operators
on $\mathfrak{h}$.

Let $C$ be a self-adjoint positive operator in $\mathfrak{h}$.
Then, for any $x,y\in\mathcal{D}( C) $ we set $\langle x,y\rangle
_{C}=\langle x,y\rangle+\langle Cx,Cy\rangle$ and $ \Vert x\Vert
_{C}=\sqrt{\langle x,x\rangle_{C}}$.
As usual, $L^{2}( \mathbb{P},\mathfrak{h}) $ stands for the set of all
square integrable random variables from $( \Omega,\mathfrak{F},\mathbb
{P})$ to $ ( \mathfrak{h},\mathfrak{B}( \mathfrak{h}) )$.
We write $L_{C}^{2}( \mathbb{P},\mathfrak{h}) $ for the set of all $\xi
\in L^{2}( \mathbb{P},\mathfrak{h}) $ satisfying $\xi\in\mathcal{D}(
C) $ a.s. and $\mathbb{E} \Vert\xi\Vert_{C}^{2}<\infty$.
The function $\pi_{C}\dvtx\mathfrak{h\rightarrow h}$ is defined by
$\pi_{C}( x) = x$ if $x\in\mathcal{D}( C)$ and $\pi_{C}( x) = 0$
whenever $x\notin\mathcal{D}( C)$.
In the sequel, the letter $K$ denotes generic constants.

\section{Adjoint quantum master equation}
\label{secAdQME}

We begin by presenting in detail the notion of $C$-solution to (\ref{2}).
\begin{condition}
\label{HipN3}
Suppose that $C$ is a self-adjoint positive operator in $\mathfrak{h}$
such that
$\mathcal{D}( C)$ is a subset of the domains of
$G, L_{1}, L_{2}, \ldots,$ and
the maps $G\circ\pi_{C}, L_{1}\circ\pi_{C}, L_{2}\circ\pi_{C},
\ldots$ are measurable.
\end{condition}

\begin{definition}
\label{definicion2}
Let Hypothesis~\ref{HipN3} hold. Assume that $\mathbb{I}$ is either $[
0,\infty[ $ or $[ 0,T] $, with $T\in\mathbb{R}_{+}$.
An $\mathfrak{h}$-valued adapted process $( X_{t}( \xi) ) _{t \in
\mathbb{I}}$ with continuous sample paths
is
called strong $C$-solution of (\ref{2}) on $\mathbb{I}$ with initial
datum $\xi$
if and only if
for all $t\in\mathbb{I}$:
\begin{itemize}
\item$\mathbb{E}\Vert X_{t}( \xi) \Vert^{2}\leq\mathbb{E}\Vert\xi
\Vert^{2}$, $X_{t}( \xi) \in\mathcal{D}( C) $ a.s.,
$
\sup_{s\in[ 0,t] }\mathbb{E}\Vert C X_{s}( \xi) \Vert^{2} < \infty.
$

\item
$
X_{t}( \xi) =\xi+\int_{0}^{t}G\pi_{C}( X_{s}( \xi
) ) \,ds+\sum_{k=1}^{\infty}\int_{0}^{t}L_{k}\pi_{C}(
X_{s}( \xi) ) \,dW_{s}^{k}$
$\mathbb{P}$-a.s.

\end{itemize}
\end{definition}

\begin{notation}
\label{notX}
The symbol $X ( \xi)$ will be reserved for the strong $C$-solution of
(\ref{2}) with initial datum $\xi$.
\end{notation}

\begin{remark}
\label{notaMedibilidad}
Suppose that $C$ is a self-adjoint positive operator in $\mathfrak{h}$,
together with
$A \in\mathfrak{L}( ( \mathcal{D}( C), \mbox{$\Vert\cdot\Vert_{C}$})
,\mathfrak{h})$.
Then $A \circ\pi_{C}\dvtx \mathfrak{h} \rightarrow\mathfrak{h}$ is
measurable whenever $\mathfrak{h}$ is equipped with its Borel $\sigma
$-algebra (see, e.g.,~\cite{FagnolaMora2010} for details).
\end{remark}

We now make more precise our basic assumptions, that is hypothesis (H).

\begin{condition}
\label{HipN5}
Suppose that Hypothesis~\ref{HipN3} holds.
In addition, assume:
\begin{longlist}[(H2.3)]
\item[(H2.1)] The operator $G$ belongs to $\mathfrak{L}( ( \mathcal{D}(
C), \mbox{$\Vert\cdot\Vert_{C}$}) ,\mathfrak{h})$.
\item[(H2.2)] For all $x\in\mathcal{D}( C) $,
$
2\Re\langle x,Gx\rangle+{\sum_{k=1}^{\infty}}\Vert L_{k}x\Vert^{2} = 0
$.

\item[(H2.3)]
Let $\xi\in L_{C}^{2}( \mathbb{P}, \mathfrak{h} )$ be $\mathfrak
{F}_{0}$-measurable.
Then for all $T > 0$, (\ref{2}) has a unique strong $C$-solution on $[
0,T]$ with initial datum $\xi$.
\end{longlist}
\end{condition}

\begin{remark}
\label{nota1}
Let $A$ be a closable operator in $\mathfrak{h}$
whose domain is contained in $\mathcal{D}( C)$,
where $C$ is a
self-adjoint positive operator in $\mathfrak{h}$. Applying the closed
graph theorem we obtain $A\in\mathfrak{L}( ( \mathcal{D}( C) ,\mbox{$\Vert
\cdot\Vert_{C}) ,\mathfrak{h}$}) $,
which leads to a sufficient condition for (H2.1).
\end{remark}

\begin{remark}
\label{nota3}
Let $C$ be a self-adjoint positive operator in $\mathfrak{h}$ such that
$\mathcal{D}( C) \subset\mathcal{D}( G)$.
Assume that
$
2\Re\langle x, Gx\rangle+\sum_{k=1}^{\infty}\Vert L_{k}x \Vert
^{2}\leq0
$
for all $x \in\mathcal{D}( G)$.
Then the numerical range of $G$ is contained in the left half-plane of
$\mathbb{C}$, and so $G$ is closable.
Therefore $G \in\mathfrak{L}( ( \mathcal{D}( C), \Vert\cdot\Vert
_{C}) ,\mathfrak{h})$
by Remark~\ref{nota1}.
\end{remark}

Using arguments given in Section~\ref{subsecuniqueness}
we prove the following theorem,
establishing the uniqueness of the solution of (\ref{41}).

\begin{definition}
\label{definicion3}
Suppose that
$A \in\mathfrak{L}( \mathfrak{h})$
and that
$C$ is a self-adjoint positive operator in $\mathfrak{h}$.
A family of operators $( \mathcal{A}_{t} )_{t \geq0}$ belonging to
$\mathfrak{L}( \mathfrak{h})$ is a $C$-solution of (\ref{41}) with
initial datum $A$ iff
$\mathcal{A}_{0} = A$
and
for all $t \geq0$:
\begin{longlist}[(b)]
\item[(a)]
$
\frac{d}{dt} \langle x, \mathcal{A}_{t} y \rangle= \langle x, \mathcal
{A}_{t} G y \rangle+ \langle G x, \mathcal{A}_{t} y \rangle+ \sum
_{k=1}^{\infty} \langle L_{k} x, \mathcal{A}_{t} L_{k} y \rangle
$
for all $x,y \in\mathcal{D}(C )$.

\item[(b)]
${\sup_{s \in[ 0,t]}} \Vert\mathcal{A}_{s} \Vert_{\mathfrak{L}(
\mathfrak{h})} < \infty$.
\end{longlist}
\end{definition}
\begin{theorem}
\label{teorema3}
Suppose that Hypothesis~\ref{HipN5} holds.
Let $A$ belong to $\mathfrak{L}( \mathfrak{h})$.
Then, for every nonnegative real number $t$ there exists a unique
$\mathcal{T}_{t}( A ) $ in $\mathfrak{L}( \mathfrak{h})$ such that for
all $x,y$ in $\mathcal{D}(C )$,
%
%
\begin{equation}
\label{42}
\langle x, \mathcal{T}_{t}( A ) y \rangle= \mathbb{E} \langle X_{t} (
x ), A X_{t} ( y ) \rangle.
\end{equation}
Moreover,
any $C$-solution of (\ref{41}) with initial datum $A$ coincides with
$ \mathcal{T} ( A ) $,
and
$
\Vert\mathcal{T}_{t}( A ) \Vert_{\mathfrak{L}( \mathfrak{h})}
\leq
\Vert A \Vert_{\mathfrak{L}( \mathfrak{h})}$
for all $t \geq0$.
\end{theorem}
\begin{pf}
The proofs fall naturally into Lemmata~\ref{lema41} and~\ref{lema42}.
\end{pf}

As a by-product of our proof of
the existence of solutions to (\ref{3}),
we ``construct'' a solution to (\ref{41}),
and so
Theorem~\ref{teorema3}
leads to Theorem~\ref{teorema10}.
\begin{theorem}
\label{teorema10}
Let Hypothesis~\ref{HipN5} hold.
Suppose that $A \in\mathfrak{L}( \mathfrak{h})$ and
that $ \mathcal{T}_{t}( A )$ is as in Theorem~\ref{teorema3}.
Then
$( \mathcal{T}_{t}( A ) )_{t \geq0}$ is the unique $C$-solution of
(\ref{41}) with initial datum $A$.
\end{theorem}
\begin{pf}
Lemmata~\ref{lema26} and~\ref{lema20} shows that
$( \mathcal{T}_{t}( A ) )_{t \geq0}$ is a $C$-solution of (\ref{41})
with initial datum $A$.
Theorem~\ref{teorema3} now completes the proof.
\end{pf}

\begin{remark}
In~\cite{MoraJFA2008},
C. M. Mora developed the existence and uniqueness of the solution
to (\ref{41}) with $A$ unbounded,
as well as its probabilistic representation.
Thus taking $A \in\mathfrak{L}( \mathfrak{h})$,
Corollary 14 of~\cite{MoraJFA2008} established
the statement of Theorem~\ref{teorema10}
under assumptions including the existence of
an orthonormal basis $( e_{n} )_{n \in\mathbb{N}}$ of $\mathfrak{h}$
that satisfies, for example,
$G e_n, L_k e_n \in\mathcal{D}( C)$
and
${\sup_{n\in\mathbb{Z}_{+}}}\Vert CP_{n}x\Vert
\leq\Vert C x\Vert$ for all $x \in\mathcal{D}( C)$, where $P_{n}$ is
the orthogonal projection
of $\mathfrak{h}$ over the linear manifold spanned by
$e_{0},\ldots, e_{n}$.
In Theorem~\ref{teorema10} we remove this basis,
extending the range of applications.
\end{remark}

\begin{remark}
\label{nota7}
Suppose that
$
2\Re\langle x,Gx\rangle+{\sum_{k=1}^{\infty}}\Vert L_{k}x\Vert^{2}
\leq0
$
for all $x \in\mathcal{D}(G )$.
Let
$G$ be the infinitesimal generator of a $C_0$-semigroup of contractions.
Define the sequence $( \mathcal{T}^{(n)} )_{n \geq0}$
of linear contractions on $\mathfrak{L}( \mathfrak{h})$
by
\[
\bigl\langle u, \mathcal{T}_{t}^{(n+1)}( A ) v \bigr\rangle
=
\langle e^{G t } u, A e^{G t } v \rangle
+
\sum_{k=1}^{\infty} \int_{0}^{t} \bigl\langle L_{k} e^{G ( t -s )} u,
\mathcal{T}_{s}^{(n)} ( A ) L_{k} e^{G ( t -s )} v \bigr\rangle
\,ds,
\]
where $u,v \in\mathcal{D}(G ) $, $A\in\mathfrak{L}( \mathfrak{h})$,
and $\mathcal{T}^{(-1)} = 0$.
A. M. Chebotarev proved that Picard's successive approximations $\mathcal
{T}^{(n)}$
converge as $n \rightarrow\infty$
to a quantum dynamical semigroup $\mathcal{T}^{(\min)}$ which is a weak
solution to (\ref{41}); see, for example,~\cite{Chebotarev2000,Fagnola1999}.
Holevo~\cite{Holevo1996} developed the probabilistic representation of
$\mathcal{T}^{(\min)}$
under restrictions, including that
$G$ and $G^{\ast}$ are the infinitesimal generators of $C_0$-semigroup
of contractions.
From Chebotarev and Fagnola~\cite{ChebFagn1998} we have that
$\mathcal{T}^{(\min)}_{t}( I ) = I$ for any $t \geq0$,
provided that there exists a self-adjoint positive operator $C$ in
$\mathfrak{h}$
and a linear manifold $\mathfrak{D} \subset\mathcal{D}(G )$ which is
a core for $C$ such that:
(i) The semigroup generated by $G$ leaves invariant $\mathfrak{D}$;
and
(ii) For some $\gamma>0$,
$
2\Re\langle C^{2} x, Gx\rangle+\sum_{k=1}^{\infty}\Vert C L_{k}x
\Vert^{2}\leq\alpha\Vert x\Vert_{C}^{2}
$
for all $x \in( \gamma I - G )^{-1} ( \mathfrak{D})$; see also \cite
{ChebGarQue1998,Fagnola1999}.
This implies the uniqueness (in the semigroup sense) of the solution
to (\ref{41}) with $A$ bounded; see, for example,~\cite{Fagnola1999}.

In addition to its proof,
the main novelty of Theorem~\ref{teorema10} is that
we do not assume properties like
$G$ are the infinitesimal generators of a semigroup
and
condition (i),
which involves the study of invariant sets for $\exp(Gt )$.
The latter is not an easy problem in general.
\end{remark}


\section{Probabilistic representations of regular density operators}
\label{secprob-rep}


The following notion of a regular density operator was introduced by
Chebotarev, Garc\'{\i}a and Quezada~\cite{ChebGarQue1998} to investigate
the identity preserving property of minimal quantum dynamical semigroups.\vadjust{\goodbreak}

\begin{definition}
\label{def2}
Let $C$ be a self-adjoint positive operator in $\mathfrak{h}$. An
operator $\varrho$ belonging to $\mathfrak{L}_{1}^{+}( \mathfrak{h} )$
is called $C$-regular iff
$
\varrho=\sum_{n\in\mathfrak{I}}\lambda_{n}\vert u_{n}\rangle\langle
u_{n}\vert
$
for some countable set $\mathfrak{I}$, summable nonnegative real
numbers $( \lambda_{n}) _{n\in\mathfrak{I}}$ and family $( u_{n}) _{n\in
\mathfrak{I}}$ of elements of $\mathcal{D}( C) $, which together satisfy
$
\sum_{n\in\mathfrak{I}}\lambda_{n}\Vert Cu_{n}\Vert^{2}<\infty
$.
We write $\mathfrak{L}_{1,C}^{+}( \mathfrak{h}) $ for the set of all
$C$-regular density operators.
\end{definition}

We next formulate the concept of $C$-regular operators in terms of
random variables.
This characterization of $\mathfrak{L}_{1,C}^{+}( \mathfrak{h}) $
complements those given in
\cite{ChebGarQue1998} using operator theory; see also~\cite{Chebotarev2000}.
\begin{theorem}
\label{teorema4}
Suppose that $C$ is a self-adjoint positive operator in $\mathfrak{h}$.
Let $\varrho$ be a linear operator in $\mathfrak{h}$. Then $\varrho$
is $C$-regular if and only if $\varrho= \mathbb{E}\vert\xi\rangle
\langle\xi\vert$
for some $\xi\in L_{C}^{2}( \mathbb{P},\mathfrak{h}) $.
Moreover, $\mathbb{E}\vert\xi\rangle\langle\xi\vert$ can be
interpreted as a Bochner integral in both $\mathfrak{L}_{1}( \mathfrak
{h}) $ and $\mathfrak{L}( \mathfrak{h}) $.
\end{theorem}
\begin{pf}
The proof is divided into Lemmata~\ref{lema6} and~\ref{lema53}.
\end{pf}

By the following theorem,
the mean values of a large number of unbounded observables are well posed
when the density operators are $C$-regular.
Theorem~\ref{teorema8} also provides probabilistic interpretations of
these expected values.

\begin{theorem}
\label{teorema8}
Suppose that $C$ is a self-adjoint positive operator in $\mathfrak{h}$,
and fix
$\varrho=\mathbb{E}\vert\xi\rangle\langle\xi\vert$ with $\xi\in
L_{C}^{2}( \mathbb{P},\mathfrak{h}) $.
Then:
\begin{longlist}[(b)]
\item[(a)]
The range of $\varrho$ is contained in $\mathcal{D}( C)$
and $
C\varrho= \mathbb{E}\vert C \xi\rangle\langle\xi\vert
$.

\item[(b)]
Consider
$A \in\mathfrak{L}( ( \mathcal{D}( C) ,\mbox{$\Vert\cdot\Vert_{C}$})
,\mathfrak{h}) $,
and let $B$ be a densely defined linear operator in $\mathfrak{h}$ such that
$\mathcal{D}( C) \subset\mathcal{D}( B^{\ast})$.
Then
$A\varrho B$ is densely defined and bounded.
The unique bounded extension of $A\varrho B$ belongs to $\mathfrak
{L}_{1}( \mathfrak{h}) $ and is equal to $\mathbb{E}\vert A\xi\rangle
\langle B^{\ast} \xi\vert$,
where
$\mathbb{E}\vert A\xi\rangle\langle B^{\ast} \xi\vert$ is a well
defined Bochner integral
in both $\mathfrak{L}_{1}( \mathfrak{h}) $ and $\mathfrak{L}( \mathfrak
{h}) $.
Moreover,
\[
\operatorname{tr}( A\varrho B ) =\mathbb{E}\langle B^{\ast} \xi,A\xi
\rangle.
\]
\end{longlist}
\end{theorem}
\begin{pf}
Deferred to Section~\ref{subsecteorema8}.
\end{pf}

\section{Quantum master equation}
\label{secQME}

We first deduce that
(\ref{i4}) defines a density operator.
\begin{theorem}
\label{teor7}
Let Hypothesis~\ref{HipN5} hold.
Then, for every $t\geq0 $ there exists a unique operator
$\rho_{t} \in\mathfrak{L} ( \mathfrak{L}_{1}( \mathfrak{h}) ) $
such that for each $C$-regular operator $\varrho$,
%
%
\begin{equation}
\label{31}
\rho_{t}( \varrho) =\mathbb{E}\vert X_{t}( \xi)
\rangle\langle X_{t}( \xi) \vert,
\end{equation}
where $\xi$ is an arbitrary random variable in $L_{C}^{2}( \mathbb
{P},\mathfrak{h}) $ satisfying
$\varrho=\mathbb{E} \vert\xi\rangle\langle\xi\vert$.
Here
$ X( \xi) $ is the strong $C$-solution of (\ref{2}) with initial datum
$\xi$,
and we can interpret
$
\mathbb{E}\vert X_{t}( \xi) \rangle\langle X_{t}( \xi) \vert
$
as a Bochner integral in $\mathfrak{L}_{1}( \mathfrak{h}) $ as well as
in $\mathfrak{L}( \mathfrak{h}) $.
Moreover,
$
\Vert\rho_{t} \Vert_{\mathfrak{L} ( \mathfrak{L}_{1}( \mathfrak{h})
)} \leq1
$
for all $t \geq0$.
\end{theorem}
\begin{pf}
Deferred to Section~\ref{subsecteor7}.
\end{pf}

\begin{notation}
\label{notrho}
From now on, $\rho_{t}$ stands for the operator given by (\ref{31}).
\end{notation}

The next theorem says that
the expected value $\mathbb{E}$
commutes with
the action of $\rho_{t}$
on random $C$-regular pure density operators.
\begin{theorem}
\label{teor9}
Assume that Hypothesis~\ref{HipN5} holds.
Let $\varrho= \mathbb{E}\vert\xi\rangle\langle\xi\vert$,
with $\xi\in$ $L_{C}^{2}( \mathbb{P},\mathfrak{h}) $.
Then
$
\mathbb{E}\rho_{t}( \vert\xi\rangle\langle\xi\vert)
=
\rho_{t}( \varrho)
$
for all $t\geq0$.
\end{theorem}
\begin{pf}
Deferred to Section~\ref{subsecteor9}.
\end{pf}

We now summarize
some relevant properties of the family of linear operators $ ( \rho_{t}
)_{t \geq0}$.
\begin{theorem}
\label{teor8}
Adopt Hypothesis~\ref{HipN5}.
Then $( \rho_{t})_{t\geq0}$ is a semigroup of contractions such that
$\rho_{t}(
\mathfrak{L}_{1}^{+}( \mathfrak{h}) ) \subset\mathfrak{L}_{1}^{+}(
\mathfrak{h}) $,
$\rho_{t}( \mathfrak{L}_{1,C}^{+}( \mathfrak{h}) ) \subset\mathfrak
{L}_{1,C}^{+}( \mathfrak{h}) $,
and
for all $\varrho\in\mathfrak{L}_{1,C}^{+}( \mathfrak{h}) $,
%
%
\begin{equation}
\label{310}
{\lim_{s\rightarrow t}\operatorname{tr}}\vert\rho_{s}( \varrho) -\rho
_{t}( \varrho) \vert= 0 .
\end{equation}
\end{theorem}
\begin{pf}
The proof is divided into Lemmata~\ref{lema14},~\ref{lema24} and~\ref{lema25}.
\end{pf}

The analysis outlined in Section~\ref{subsecexistence}
leads to our first main theorem,
which asserts that
$\mathbb{E}\vert X_{t}( \xi)
\rangle\langle X_{t}( \xi) \vert$ satisfies (\ref{3})
in both senses, integral and $\mathfrak{L}_{1}( \mathfrak{h}) $-weak,
whenever $\varrho= \mathbb{E} \vert\xi\rangle\langle\xi\vert$ is $C$-regular.

\begin{condition}
\label{HipN1}
The operators $G, L_{1}, L_{2}, \ldots$ are closable.
\end{condition}
\begin{theorem}
\label{teor10}
Let Hypotheses~\ref{HipN5} and~\ref{HipN1} hold.
Suppose that $\varrho$ is $C$-regular.
Then for all $t\geq0$,
%
%
\begin{equation}
\label{311}
\rho_{t}( \varrho) =\varrho+\int_{0}^{t}\Biggl( G\rho_{s}(
\varrho) +\rho_{s}( \varrho) G^{\ast}+\sum_{k=1}^{\infty
}L_{k}\rho_{s}( \varrho) L_{k}^{\ast}\Biggr) \,ds,
\end{equation}
where
we understand the above integral in the sense of the Bochner integral
in $\mathfrak{L}_{1}( \mathfrak{h}) $.
Moreover, for any $A\in\mathfrak{L}( \mathfrak{h}) $ and $t\geq0$,
%
%
\begin{equation}
\label{312}
\frac{d}{dt}\operatorname{tr}( A\rho_{t}( \varrho) )
=
\operatorname{tr}\Biggl(
A\Biggl( G\rho_{t}( \varrho) +\rho_{t}( \varrho)G^{\ast}
+\sum_{k=1}^{\infty}L_{k}\rho_{t}( \varrho) L_{k}^{\ast}\Biggr)
\Biggr).
\end{equation}
\end{theorem}
\begin{pf}
Deferred to Section~\ref{subsecteor10}.
\end{pf}

\begin{remark}
Let $G, L_{1}, L_{2}, \ldots$ be densely defined.
Then Hypothesis~\ref{HipN1} is equivalent to saying that
$G^{\ast}, L_{1}^{\ast}, L_{2}^{\ast}, \ldots$ are densely defined.\vadjust{\goodbreak}
\end{remark}

The second main theorem of this paper
establishes that
under Hypothesis~\ref{HipN5},
$\rho_{t}( \varrho)$
is the unique solution of (\ref{312}) in the semigroup sense.
Its proof is based on
arguments given in Section~\ref{subsecuniqueness}.

\begin{definition}
\label{defSemigroupSol}
A semigroup $( \widehat{\rho}_{t}) _{t\geq0}$
of bounded operators on $\mathfrak{L}_{1}( \mathfrak{h}) $
is called semigroup $C$-solution of (\ref{3}) if and only if:
\begin{longlist}[(iii)]
\item[(i)] For each nonnegative real number $T$, ${\sup_{t\in[ 0,T]
}}\Vert\widehat{\rho}_{t}\Vert_{\mathfrak{L}( \mathfrak{L}_{1}(
\mathfrak{h}) ) }<\infty$.

\item[(ii)] For any $x\in\mathcal{D}( C) $ and $A\in\mathfrak{L}(
\mathfrak{h}) $,
the function $t\mapsto\operatorname{tr}( \widehat{ \rho}_{t}( \vert
x\rangle\langle x\vert) A) $ is continuous.

\item[(iii)]
$
\lim_{t\rightarrow0+}
(
\operatorname{tr}( A \widehat{\rho}_{t} ( \vert x\rangle\langle x\vert
) )
-
\operatorname{tr}( A \vert x\rangle\langle x\vert)
) / t
=
\langle x, A Gx\rangle+\langle Gx, A x\rangle
+\break
\sum_{k=1}^{\infty}\langle L_{k}x, A L_{k}x\rangle
$
whenever $x\in\mathcal{D}( C) $ and $A\in\mathfrak{L}( \mathfrak{h}) $.
\end{longlist}
\end{definition}
\begin{theorem}
\label{teorema9}
Let Hypothesis~\ref{HipN5} hold.
Then $( \rho_{t})_{t\geq0}$ is
the unique semigroup $C$-solution of (\ref{3}).
\end{theorem}
\begin{pf}
Deferred to Section~\ref{subsecteorema9}.
\end{pf}

Theorems~\ref{teor10} and~\ref{teorema9},
together with Theorem~\ref{teorema8},
show that
the mean values of the observables
with respect to the solutions of the quantum master equations
are well posed
in many physical situations.
Moreover,
Theorems~\ref{teorema8},~\ref{teor10} and~\ref{teorema9}
allow us to make rigorous some explicit computations
concerning the evolution of unbounded observables,
like the following Ehrenfest-type theorem.
\begin{theorem}
\label{teorema5}
Assume the setting of Example~\ref{exmeasurement}.
Then
$( \rho_{t})_{t\geq0}$ is
the unique semigroup $( P^2 + Q^2 )$-solution of (\ref{3}).
If
$\varrho\in\mathfrak{L}_{1, P^2 + Q^2}^{+}( L^{2}( \mathbb{R},\mathbb
{C}) ) $,
then for all $t \geq0$,
%
%
\begin{equation}
\label{31n}\quad
\frac{d}{dt}\operatorname{tr}( Q \rho_{t}( \varrho) )
=
\frac{1}{m} \operatorname{tr}( P \rho_{t}( \varrho) ) ,\qquad
\frac{d}{dt} \operatorname{tr}( P \rho_{t}( \varrho) )
=
- 2c \operatorname{tr}( Q \rho_{t}( \varrho) ) .
\end{equation}
\end{theorem}
\begin{pf}
Deferred to Section~\ref{subsecteorema5}.
\end{pf}

\begin{remark}
\label{nota8}
A novelty of this paper lies in
the use of probabilistic methods for proving
Theorems~\ref{teor10} and~\ref{teorema9}.
In order to adopt a purely Operator Theory viewpoint,
we now return to Remark~\ref{nota7}.
Let $( \mathcal{T}_{\ast t} ) _{t\geq0}$ be the semigroup on $\mathfrak
{L}_{1} ( \mathfrak{h})$
whose adjoint semigroup is $( \mathcal{T}^{(\min)}_{t} ) _{t\geq0}$;
that is, $\mathcal{T}_{\ast}$ is the predual semigroup of $\mathcal
{T}^{(\min)}$.
In case $\mathcal{T}^{(\min)}$ leaves invariant the identity operator,
the linear span of
$ \{ \vert x\rangle\langle y\vert\dvtx x,y\in\mathcal{D}( G ) \}$
is a core for the infinitesimal generator of $( \mathcal{T}_{\ast t} )
_{t\geq0}$,
which is denoted by $\mathcal{L}_{\ast}$ for simplicity of notation;
see, for example, Proposition 3.32 of~\cite{Fagnola1999}.
Then,
under conditions (i) and (ii) given in Remark~\ref{nota7},
$( \mathcal{T}_{\ast t}) _{t\geq0}$ is the unique strongly continuous
semigroup on $\mathfrak{L}_{1}( \mathfrak{h}) $ satisfying a version of
(\ref{3}) for all
$\varrho= \vert x \rangle\langle y \vert$
with $x,y \in\mathcal{D}( G) $.
In order to establish
$ \mathcal{T}_{\ast t}( \mathfrak{L}_{1,C}^{+}( \mathfrak{h}) )
\subset\mathfrak{L}_{1,C}^{+}( \mathfrak{h}) $
as well as
the assertions of Theorem~\ref{teor10},
we have to prove first that
$\mathfrak{L}_{1,C}^{+}( \mathfrak{h}) \subset\mathcal{D}( \mathcal
{L}_{\ast} ) $.
If we are able to do it,
then
$\mathfrak{L}_{1,C}^{+}( \mathfrak{h}) $ is an invariant set for $
\mathcal{T}_{\ast t} $
provided that
%
%
\begin{equation}
\label{i10}
\sup_{n \in\mathbb{N}} \bigl| \operatorname{tr}\bigl( \mathcal{T}^{(\min)} \bigl( n(n
+ C^2 ) ^{-1} C^2 \bigr) \varrho\bigr) \bigr|
<
\infty
\end{equation}
for any $\varrho\in\mathfrak{L}_{1,C}^{+}( \mathfrak{h})$.
When $C$ is invertible,
(\ref{i10}) follows from
%
%
\begin{equation}
\label{i11} \bigl\| C^{-1} \mathcal{T}^{(\min)} \bigl( n(n + C^2 ) ^{-1}C^2
\bigr) C^{-1} \bigr\| \leq K \| n(n + C^2 ) ^{-1} \| .
\end{equation}
A careful reading of~\cite{GarQue1998} reveals that
for any $A \in\mathfrak{L}( \mathfrak{h})$ we have
%
%
\begin{equation}
\label{i12}
\bigl\| C^{-1} \mathcal{T}^{(\min)} ( A ) C^{-1} \bigr\|
\leq
K \| C^{-1} A C^{-1} \|
\end{equation}
under assumptions of type (\ref{i2}),
together with
$ \exp( G t )$,
leaves invariant a core of $C$ contained in $\mathcal{D}( G )$; see
also~\cite{ChebGarQue1998,Chebotarev2000}.
This gives (\ref{i11}),
and so (\ref{i10}) holds.
Under the same assumptions,\vspace*{1pt}
an alternative is to obtain (\ref{i12}) by proving
$
\| C^{-1} \mathcal{T}^{(n)} ( A )\times C^{-1} \|
\leq
K \| C^{-1} A C^{-1} \|
$
directly from the definition of $\mathcal{T}^{(n)}$,
but with effort.
Here
$\mathcal{T}^{(n)}$ is as in Remark~\ref{nota7}.
Finally,
to establish (\ref{311}) and (\ref{312})
we have to get that
$
\mathcal{L}_{\ast}( \varrho) =G\varrho+\varrho
G^{\ast}+\sum_{k=1}^{\infty}L_{k}\varrho L_{k}^{\ast}
$
for any $\varrho\in\mathfrak{L}_{1,C}^{+}( \mathfrak{h})$.
\end{remark}

\section{Regular stationary solutions of quantum master equations}
\label{secStatSol}

This section is devoted to objective (O2).
In this direction,
the next theorem provides the representation of the density operator at
time $t$
as the average of all pure states $\vert Y_{t}\rangle\langle Y_{t}\vert$
associated to the nonlinear stochastic Schr\"{o}dinger equation (\ref{5}).
This model has a sound physical basis; see, for example,
\cite
{BarchielliBelavkin1991,BarchielliGregoratti2009,GoughSobolev2004,ScottMilburn2001}.

\begin{definition}
Let $C$ satisfy Hypothesis~\ref{HipN3}.
Suppose that $\mathbb{I}$ is either $[ 0,+\infty[ $ or $[ 0,T] $
provided $T\in[ 0,+\infty[ $.
We say that
$( \mathbb{Q}, ( Y_{t})_{t\in\mathbb{I}},( B_{t})_{t\in\mathbb{I}})$
is a $C$-solution of (\ref{5}) with initial distribution $\theta$ on
$\mathbb{I}$ if and only if:
\begin{itemize}
\item$B = ( B^{k} )_{k \in\mathbb{N}}$ is a sequence of real valued
independent Brownian motions on
the filtered complete probability space $( \Omega,\mathfrak{F},(
\mathfrak{F}_{t}) _{t\in\mathbb{I}},\mathbb{Q}) $.

\item$( Y_{t}) _{t\in\mathbb{I}}$ is an $\mathfrak{h}$-valued
process with continuous sample paths such that the law of $Y_{0}$ coincides
with $\theta$ and $\mathbb{Q}( \Vert Y_{t}\Vert=1$
for all $t\in\mathbb{I}) =1$.

\item For every $t\in\mathbb{I}\dvtx Y_{t}\in\mathcal{D}( C)$ $\mathbb{Q}$-a.s.
and $\sup_{s\in[ 0,t] }\mathbb{E}_{\mathbb{Q}}\Vert CY_{s}\Vert
^{2}<\infty$.\vspace*{1pt}

\item$\mathbb{Q}$-a.s.,
$
Y_{t}=Y_{0}+\int_{0}^{t}G( \pi_{C} (Y_{s})) \,ds+\sum_{k=1}^{\infty}\int
_{0}^{t}L_{k}( \pi_{C} ( Y_{s} ) ) \,dB_{s}^{k}
$
for all $t\in\mathbb{I}$.
\end{itemize}
\end{definition}
\begin{theorem}
\label{corolario2}
Suppose that Hypothesis~\ref{HipN5} holds.
Let
$
\varrho=\int_{\mathfrak{h}} \vert y\rangle\langle y\vert\theta( dy)
$,
with $\theta$ probability measure over $\mathfrak{h}$ satisfying
$
\theta( \mathcal{D}( C) \cap\{ x\in\mathfrak{h}\dvtx\Vert
x\Vert=1\} ) =1
$
and
$
\int_{\mathfrak{h}}\Vert Cx\Vert^{2}\theta( dx) <\infty
$.
Then for all $t \geq0$,
\[
\rho_{t}( \varrho) =\mathbb{E}_{ \mathbb{Q}} \vert Y_{t}\rangle\langle
Y_{t}\vert,
\]
where
$\rho_{t}( \varrho)$ is defined by (\ref{31}),
and
$( \mathbb{Q},( Y_{t}) _{t\geq0},( B_{t}) _{t\geq0}) $
is the $C$-solution of (\ref{5}) with initial law $\theta$.
\end{theorem}
\begin{pf}
Deferred to Section~\ref{subseccorolario2}.
\end{pf}

\begin{remark}
\label{nota6}
Let $\theta$ be as in Theorem~\ref{corolario2}.
Suppose that Hypothesis~\ref{HipN5} holds.
Then,
we can use the same arguments as in the proof of Theorem 1 of \cite
{MoraReAAP2008}
for establishing that (\ref{5}) has a unique (in the probabilistic
sense) $C$-solution
$(\mathbb{Q},( Y_{t}) _{t\geq0},( B_{t})_{t\geq0})$
with initial law $\theta$.
\end{remark}

\begin{remark}
From Theorems~\ref{teorema8} and~\ref{corolario2}
we obtain that
the expected value of
$A\in\mathfrak{L}( ( \mathcal{D}( C), \mbox{$\Vert\cdot\Vert_{C}$})
,\mathfrak{h})$
at time $t$ is equal to
$
\mathbb{E} \langle Y_{t} , A Y_{t} \rangle
$.
This gives theoretical support for
the numerical computation of the mean value of an observable $A$ at
time $t$ through (\ref{5}),
which is the principal method for computing efficiently
$
\operatorname{tr}( A\rho_{t}( \varrho) )
$; see, for example,~\cite{BreuerPetruccione2002,MoraAAP2005,Percival1998}.
\end{remark}

Our third main theorem deals with the existence of regular stationary
states for~(\ref{3}).
This is a step forward in the understanding of the long-time behavior
of unbounded observables.

\begin{condition}
\label{Hip2}
Let Hypothesis~\ref{HipN5} hold.
Assume
the existence of a probability measure $\Gamma$ on
$\mathfrak{B}( \mathfrak{h}) $ such that:
$\Gamma( \operatorname{Dom}( C) \cap\{ x\in\mathfrak{h}\dvtx\Vert x\Vert=1\} ) =1$,
$\int_{\mathfrak{h}} \Vert C z\Vert^{2}\Gamma( dz) <\infty$
and
%
%
\begin{equation}
\label{IM1}%
\Gamma( A) =\int_{\mathfrak{h}}P_{t}( x,A)
\Gamma( dx)
\end{equation}
for any $t\geq0$ and $A\in\mathfrak{B}( \mathfrak{h}) $.
Here
$
P_{t}( x,A) = \mathbb{Q}_{x}( Y_{t}^{x}\in A)
$
if
$x\in \operatorname{Dom}( C)$
and
$
P_{t}( x,A) = \delta_{x}( A)
$
otherwise;
the $C$-solution of (\ref{5}) with initial data $x\in \operatorname{Dom}( C) $
is denoted by
$( \mathbb{Q}_{x},( Y_{t}^{x}) _{t\geq0},(B_{t}^{\cdot,x}) _{t\geq0}) $.
\end{condition}
\begin{theorem}
\label{teorema7}
Under Hypothesis~\ref{Hip2},
there exists a $C$-regular operator
$\varrho_{\infty}$ such that
$
\rho_{t}( \varrho_{\infty}) =\varrho_{\infty}
$
for all $t\geq0$.
\end{theorem}
\begin{pf}
Deferred to Section~\ref{subsecteorema7}.
\end{pf}

\begin{remark}
Combining the results of Section~\ref{secQME} with Theorem \ref
{teorema7} yields
the existence of a $C$-regular stationary solution to (\ref{3}).
\end{remark}

\section{Quantum oscillator}
\label{secoscillator}

In this section
we illustrate our general results with the following quantum oscillator.

\begin{example}
\label{Ejemplo1}
Consider $\mathfrak{h}=l^{2}( \mathbb{Z}_{+}) $,
together with its canonical orthonormal basis $( e_{n}) _{n\in\mathbb{Z}_{+}}$.
The closed operators $a^{\dagger}$, $a$ are given by:
for all $n \in\mathbb{Z}_{+}$ $a^{\dagger}e_{n}= \sqrt{ n+1}e_{n+1}$,
$ae_{0} = 0$
and
$ae_{n} = \sqrt{n}e_{n-1}$
if $n \in\mathbb{N}$.
Define $N=a^{\dagger}a$.


Choose
$
H=i\beta_{1}( a^{\dagger}-a) +\beta_{2}N+\beta_{3}( a^{\dagger}) ^{2}a^{2}
$
with $\beta_{1}, \beta_{2}, \beta_{3} \in\mathbb{R}$.
Let
$L_{1}=\alpha_{1}a$, $L_{2}=\alpha_{2}a^{\dagger}$, $L_{3}=\alpha_{3}N$,
$L_{4}=\alpha_{4}a^{2}$, $L_{5}=\alpha_{5}( a^{\dagger}) ^{2}$ and
$L_{6}=\alpha_{6}N^{2}$,
where $\alpha_{1},\ldots,\alpha_{6} \in\mathbb{C}$.
Set $L_{k}=0$ for any $k\geq7$,
and so take
$
G=-iH- \sum_{k=1}^{6}L_{k}^{\ast}L_{k} /2
$.
\end{example}

Example~\ref{Ejemplo1}
describes a laser-driven quantum oscillator in a Kerr medium
that interacts with a thermal bath.
In addition,
Example~\ref{Ejemplo1} unifies concrete physical systems such as
the following two basic models:
\begin{itemize}
\item
A mode with natural frequency $\omega$ of a electromagnetic field
inside of a cavity
is described by
$\beta_{2}=\omega$, $\alpha_{1}=\sqrt{A(\nu+1)}$, $\alpha_{2}=\sqrt
{A\nu}$
and
$\beta_{1}=\beta_{3}=\alpha_{k }=0$, with $k =3, \ldots, 6$.
Here,
the mode is damped with rate $\alpha_{1}$ by a thermal reservoir,
and
$\nu$ is a parametrization of the bath temperature;
see, for example, \cite
{BreuerPetruccione2002,GardinerZoller2004,WisemanMilburn2010}.

\item
A simple two-photon absorption and emission process is modeled by
$\beta_{3}\in\mathbb{R}$, $\alpha_{4}>0$, $\alpha_{5} \geq0$
and
$\beta_{1} = \beta_{2} = \alpha_{1} = \alpha_{2} =\alpha_{3} = \alpha
_{6} =0$;
see, for example,~\cite{CarboneFagGaQue2008,FagQue2005} and references therein.
\end{itemize}

The next theorem characterizes
the well-posedness of
the mean values of observables formed by a finite composition of
$a^{\dagger}$ and $a$
in transient and stationary regimes.
Important examples of such observables are
$ Q = i( a^{\dagger}+a) /\sqrt{2}$,
$P = i( a^{\dagger}-a) /\sqrt{2}$
and $N$.
%
\begin{theorem}
\label{teor14}
Assume the setting of Example~\ref{Ejemplo1},
and let $\rho_{t}( \varrho) $ be as in Theorem~\ref{teor7}.
Suppose that
$p$ is a natural number greater than or equal to~$4$.\looseness=-1

\begin{longlist}[(ii)]
\item[(i)]
Let $\vert\alpha_{4}\vert\geq\vert\alpha_{5}\vert$
and let
$\varrho\in\mathfrak{L}_{1,N^{p}}^{+}( l^{2}( \mathbb{Z}_{+}) ) $.
Then $\rho_{t}( \varrho) $ is a $N^{p}$-regular operator that
satisfies both (\ref{311}) and (\ref{312}).
Moreover, $( \rho_{t}) _{t\geq0}$ is the unique
semigroup $N^{p}$-solution of (\ref{3}).

\item[(ii)]
Suppose that
either $\vert\alpha_{4}\vert>\vert\alpha_{5}\vert$
or $\vert\alpha_{4}\vert=\vert\alpha_{5}\vert$ with
$\vert\alpha_{2}\vert^{2}-\vert\alpha_{1}\vert
^{2}+4( 2p+1) \vert\alpha_{4}\vert^{2}<0$.
Then,
there exists a $N^{p}$-regular
operator $\varrho_{\infty}$ such that
$
\rho_{t}( \varrho_{\infty}) =\varrho_{\infty}
$
for any $t\geq0$.
\end{longlist}
\end{theorem}
\begin{pf}
Deferred to Section~\ref{subsecteor14}.
\end{pf}

\begin{remark}
\label{notaSuffCond}
In the proof of Theorem~\ref{teor14} we use the following sufficient
condition for condition (H2.3),
which is developed in~\cite{FagnolaMora2010}.
\end{remark}
%
\begin{condition}
\label{HipN4}
Suppose that $C$ is a self-adjoint positive operator in
$\mathfrak{h}$
such that
$G, L_{1}, L_{2}, \ldots$
belong to $\mathfrak{L}( ( \mathcal{D}( C), \mbox{$\Vert\cdot\Vert_{C}$})
,\mathfrak{h})$,
and
$
2\Re\langle x, Gx\rangle+{\sum_{k=1}^{\infty}}\Vert L_{k}x \Vert
^{2}\leq0
$
for any $x$ in a core of $C$.
In addition, assume that for any $x$ belonging to a core of $C^{2}$,
$
2\Re\langle C^{2} x, Gx\rangle+{\sum_{k=1}^{\infty}}\Vert C L_{k}x
\Vert^{2}
\leq
K ( \Vert x\Vert_{C}^{2} + 1 )
$.
\end{condition}
%

\section{Proofs}
\label{secproofs}
%
%

\subsection{\texorpdfstring{Proof of Theorem \protect\ref{teorema3}}{Proof of Theorem 2.1}}
\label{subsecteorema3}

We first prove that (\ref{42}) defines implicitly
a bounded operator $\mathcal{T}_{t}( A ) $.

\begin{lemma}
\label{lema41}
Adopt the assumptions of Hypothesis~\ref{HipN5} with the exception of
condition \textup{(H2.2)}.
Consider $A \in\mathfrak{L}( \mathfrak{h})$.
Then for every $t \geq0$ there exists a unique $\mathcal{T}_{t}( A ) $
belonging to $\mathfrak{L}( \mathfrak{h})$ for which (\ref{42}) holds
for all $x,y$ in $\mathcal{D}(C )$.
Moreover, $\Vert\mathcal{T}_{t}( A ) \Vert\leq\Vert A \Vert$ for
any $t \geq0$.\vadjust{\goodbreak}
\end{lemma}

\begin{pf}
By Definition~\ref{definicion2},
$
\vert\mathbb{E} \langle X_{t} ( x ), A X_{t} ( y ) \rangle\vert
\leq
\Vert A \Vert\Vert x \Vert\Vert y \Vert
$
for all $x,y \in\mathcal{D}( C )$.
Hence the sesquilinear form over $\mathcal{D}( C ) \times\mathcal{D}(
C )$ given by
$
(x, y ) \mapsto\mathbb{E} \langle X_{t} ( x ), A X_{t} ( y ) \rangle
$
can be extended uniquely to a sesquilinear form $\lbrack\cdot, \cdot
\rbrack$ over $\mathfrak{h} \times\mathfrak{h}$
with the property that
$
\vert\lbrack x, y \rbrack\vert
\leq
\Vert A \Vert\Vert x \Vert\Vert y \Vert
$
for any $x,y \in\mathfrak{h}$.
There exists a unique bounded operator $\mathcal{T}_{t}( A )$ on
$\mathfrak{h}$ such that $ \vert\lbrack x, y \rbrack\vert= \langle
x, \mathcal{T}_{t}( A ) y \rangle$ for all $x,y$ in $\mathfrak{h}$.
Furthermore,
$\Vert\mathcal{T}_{t}( A ) \Vert\leq\Vert A \Vert$.
\end{pf}

Using arguments given in Section~\ref{subsecuniqueness}
we next establish the uniqueness of solutions for the adjoint quantum
master equations.

\begin{lemma}
\label{lema42}
Let Hypothesis~\ref{HipN5} hold. Assume that $( \mathcal{A}_{t} )_{t
\geq0}$ is a $C$-solution of (\ref{41}) with initial datum $A \in
\mathfrak{L}( \mathfrak{h})$. Then $\mathcal{A}_{t} = \mathcal{T}_{t}(
A )$ for all $t \geq0$, where $\mathcal{T}_{t}( A )$ is as in Therorem
\ref{teorema3}.
\end{lemma}

\begin{pf}
Using It\^o's formula
we will prove that
for all $x,y \in\mathcal{D}( C )$,
%
%
\begin{equation}
\label{45}
\mathbb{E} \langle X_{t}( x ), A X_{t} ( y ) \rangle
=
\langle x, \mathcal{A}_{t} y \rangle.
\end{equation}
This, together with Lemma~\ref{lema41}, implies
$\mathcal{A}_{t} = \mathcal{T}_{t}( A )$.

Motivated by the fact that $\mathcal{A}_t$ is only a weak solution,
we fix an orthonormal basis $( e_{n} )_{n \in\mathbb{N}}$ of $\mathfrak{h}$
and consider the function
$
F_{n} \dvtx [0 , t ] \times\mathfrak{h} \times\mathfrak{h}
\rightarrow
\mathbb{C}
$
defined by
\[
F_{n} ( s, u, v ) = \langle R_{n} \overline{u}, \mathcal{A}_{t-s}
R_{n} v \rangle,
\]
where
$R_{n} = n (n+C )^{-1}$
and
$
\bar{u}
=
\sum_{n \in\mathbb{N}} \overline{\langle e_{n}, u \rangle} e_{n}
$.
Since the range of $R_{n}$ is contained in $\mathcal{D}( C )$,
condition (a) of Definition~\ref{definicion3} yield
%
%
\begin{equation}
\label{46}
\frac{d}{ds} F_{n} ( s, u, v ) = - g ( s, R_{n} \overline{u}, R_{n} v )
\end{equation}
with
$
g ( s, x, y )
=
\langle x, \mathcal{A}_{t-s} G y \rangle
+ \langle G x, \mathcal{A}_{t-s} y \rangle
+ \sum_{k=1}^{\infty} \langle L_{k} x, \mathcal{A}_{t-s} L_{k} y
\rangle
$.
According to conditions (a), (b) of Definition~\ref{definicion3},
we have that
$t \longmapsto\langle u, \mathcal{A}_{t} v \rangle$ is continuous
for all $u,v \in\mathfrak{h}$,
and so
combining $C R_{n} \in\mathfrak{L}( \mathfrak{h})$ with Hypothesis~\ref{HipN5}
we get the uniform continuity of
$
(s, u, v )
\longmapsto
g ( s, R_{n} \overline{u}, R_{n} v )
$
on bounded subsets of $[ 0, t ] \times\mathfrak{h} \times\mathfrak{h}$.
Therefore
we can apply It\^o's formula to
$ F_{n} ( s \wedge\tau_{j}, \overline{X_{s}^{\tau_{j}}( x )},
X_{s}^{\tau_{j}}( y ) ) $,
with
$
\tau_{j} =
\inf{ \{ t \geq0\dvtx \Vert X_{t}( x ) \Vert+ \Vert X_{t}( y ) \Vert> j
\} }
$.

Fix $x,y \in\mathcal{D}( C )$.
Combining It\^o's formula with (\ref{46})
we deduce that
\[
F_{n} \bigl( t \wedge\tau_{j}, \overline{X_{t}^{\tau_{j}}( x )}, X_{t}^{\tau
_{j}}( y ) \bigr)
=
F_{n} ( 0, \overline{X_{0}( x )}, X_{0}( y ) )
+ I_{t \wedge\tau_{j}}^n + M_t .
\]
Here\vspace*{1pt} for $s \in[ 0, t ]$:
$
M_{s}
=
\sum_{k = 1}^{\infty} \int_{0}^{s \wedge\tau_{j}} \langle
R_{n} X_{r} ^{\tau_{j}} ( x ),
\mathcal{A}_{t-r}
R_{n} L_{k} X_{r}^{\tau_{j}} ( y )
\rangle \,dW^{k}_{r}
+\break
\sum_{k = 1}^{\infty} \int_{0}^{s \wedge\tau_{j}} \langle
R_{n} L_k X_{r} ^{\tau_{j}} ( x ),
\mathcal{A}_{t-r}
R_{n} X_{r}^{\tau_{j}} ( y )
\rangle \,dW^{k}_{r}
$
and
\[
I_s^n
=
\int_{0}^{s}
\bigl(
- g ( r, R_n X_{r} ( x ) , R_n X_{r} ( y ) )
+
g_{n} ( r, X_{r} ( x ) , X_{r} ( y ) )
\bigr) \,dr,
\]
the function
$
g_{n} ( r, u, v )
$
is equal to
$
\langle R_{n} u, \mathcal{A}_{t-r} R_{n} G v \rangle
+ \langle R_{n} G u, \mathcal{A}_{t-r} R_{n} v \rangle
+ \sum_{k=1}^{\infty} \langle R_{n} L_{k} u, \mathcal{A}_{t-r} R_{n}
L_{k} v \rangle$.\vadjust{\goodbreak}

We next establish the martingale property of $M_s$.
For all $r \in[0, t ]$ we have
\[
\Vert R_{n} X_{r} ^{\tau_{j}} ( x ) \Vert^2
\Vert\mathcal{A}_{t-r} \Vert^2
\Vert R_{n} L_{k} X_{r}^{\tau_{j}} ( y ) \Vert^2
\leq
j^2 \sup_{s \in[0, t ]} \Vert\mathcal{A}_{s} \Vert^2
\Vert L_{k} X_{r}^{\tau_{j}} ( y ) \Vert^2 .
\]
By (H2.1) and (H2.2),
$
\mathbb{E}
\int_{0}^{t \wedge\tau_{j}}
\sum_{k = 1}^{\infty}
|
\langle
R_{n} X_{r} ^{\tau_{j}} ( x ),
\mathcal{A}_{t-r}
R_{n} L_{k} X_{r}^{\tau_{j}} ( y )
\rangle
|^{2} \,ds
<
\infty
$.
Thus
$
(
\sum_{k = 1}^{\infty}
\int_{0}^{s \wedge\tau_{j}} \langle
R_{n} X_{r} ^{\tau_{j}} ( x ),
\mathcal{A}_{t-r}
R_{n} L_{k} X_{r}^{\tau_{j}} ( y )
\rangle \,dW^{k}_{r}
)_{s \in[ 0, t ]}
$
is a martingale.
The same conclusion can be drawn for
\[
\sum_{k = 1}^{\infty}
\int_{0}^{s \wedge\tau_{j}} \langle
R_{n} L_k X_{r} ^{\tau_{j}} ( x ),
\mathcal{A}_{t-r}
R_{n} X_{r}^{\tau_{j}} ( y )
\rangle \,dW^{k}_{r},
\]
and so $( M_s )_{s \in[ 0, t ]}$ is a martingale.
Hence
%
%
\begin{equation}
\label{415}
\mathbb{E}
\langle R_{n} X_{ t }^{\tau_{j}}( x ) ,
\mathcal{A}_{t - t \wedge\tau_{j}}
R_{n} X_{t}^{\tau_{j}}( y )
\rangle
=
\langle R_{n} x , \mathcal{ A}_{t} R_{n} y \rangle
+
\mathbb{E} I_{t \wedge\tau_{j}}^n .
\end{equation}

We will take the limit as $j \rightarrow\infty$ in (\ref{415}).
Since
$
\mathbb{E} ( {\sup_{s \in[0, t ]} }\Vert X_s( \xi) \Vert^{2} )
< \infty
$
for $\xi= x, y$
(see, e.g., Theorem 4.2.5 of~\cite{Prevot2007}),
using the dominated convergence theorem,
together with the continuity of
$t \longmapsto\langle u, \mathcal{A}_{t} v \rangle$,
we get
\[
\mathbb{E} \langle R_{n} X_{t}^{\tau_{j}}( x ), \mathcal{A}_{t- t
\wedge\tau_{j}} R_{n} X_{t}^{\tau_{j}}( y ) \rangle
\longrightarrow_{j \rightarrow\infty}
\mathbb{E} \langle R_{n} X_{t}( x ), A R_{n} X_{t} ( y ) \rangle.
\]
Applying again the dominated convergence theorem yields
$
\mathbb{E} I_{t \wedge\tau_{j}}^n
\longrightarrow_{j \rightarrow\infty}
\mathbb{E} I_{t }^n
$,
and hence letting $j \rightarrow\infty$ in (\ref{415}) we deduce that
%
%
\begin{eqnarray}
\label{414}
&& \mathbb{E} \langle R_{n} X_{t} ( x ), A R_{n} X_{t} ( y ) \rangle
- \langle R_{n} x , \mathcal{ A}_{t} R_{n} y \rangle
\nonumber\\[-8pt]\\[-8pt]
&&\qquad=
\mathbb{E} \int_{0}^{t}
\bigl(
- g ( s, R_n X_{s} ( x ) , R_n X_{s} ( y ) )
+ g_{n} ( s, X_{s} ( x ), X_{s} ( y ) )
\bigr) \,ds.\nonumber
\end{eqnarray}

Finally,
we take the limit as $n \rightarrow\infty$ in (\ref{414}).
Since
$\Vert R_{n} \Vert\leq1$
and
$ R_{n}$ tends pointwise to $I$ as $n \rightarrow\infty$,
the dominated convergence theorem yields
\[
\lim_{n \rightarrow\infty} \mathbb{E} \int_{0}^{t}
g_{n} ( s, X_{s} ( x ), X_{s} ( y ) ) \,ds
=
\mathbb{E} \int_{0}^{t}
g ( s, X_{s} ( x ), X_{s} ( y ) ) \,ds .
\]
For any $x \in\mathcal{D}( C )$,
$ \lim_{n \rightarrow\infty} C R_{n} x = C x$.
By
$\Vert C R_{n} x \Vert\leq\Vert C x \Vert$,
using the dominated convergence theorem gives
\[
\lim_{n \rightarrow\infty} \mathbb{E} \int_{0}^{t}
g ( s, R_n X_{s} ( x ) , R_n X_{s} ( y ) )
\,ds
=
\mathbb{E} \int_{0}^{t}
g ( s, X_{s} ( x ), X_{s} ( y ) ) \,ds.
\]
Thus,
letting $n \rightarrow\infty$ in (\ref{414}) we obtain (\ref{45}).
\end{pf}

\subsection{\texorpdfstring{Proof of Theorem \protect\ref{teorema8}}{Proof of Theorem 3.2}}
\label{subsecteorema8}

We begin by examining the properties of the
Bochner integral
$\mathbb{E}\vert\xi\rangle\langle\chi\vert$
when $\xi, \chi\in L^{2}( \mathbb{P},\mathfrak{h})$.

\begin{lemma}
\label{lema51}
Suppose that $\xi$ and $\chi$ belong to $L^{2}( \mathbb{P},\mathfrak
{h})$. Then $\mathbb{E}\vert\xi\rangle\langle\chi\vert$ defines an
element of $\mathfrak{L}_{1}( \mathfrak{h}) $, which moreover, is given by
%
%
\begin{equation}
\label{53}
\langle x, \mathbb{E}\vert\xi\rangle\langle\chi\vert y \rangle
=
\mathbb{E} \langle x, \xi\rangle\langle\chi, y \rangle
\end{equation}
for all $x,y \in\mathfrak{h}$.
Here, $\mathbb{E}\vert\xi\rangle\langle\chi\vert$ is well defined
as a Bochner integral with values in both
$\mathfrak{L}_{1}( \mathfrak{h}) $ and $\mathfrak{L} ( \mathfrak{h}) $.
In addition,
$
\operatorname{tr} ( \mathbb{E}\vert\xi\rangle\langle\chi\vert) =
\mathbb{E} \langle\chi, \xi\rangle$.\vadjust{\goodbreak}
\end{lemma}

\begin{pf}
We first get
$\mathbb{E}\vert\xi\rangle\langle\chi\vert
\in
\mathfrak{L}_{1}( \mathfrak{h}) $.
Since the image of $\vert\xi\rangle\langle\chi\vert$ lies in the set
of all rank-one operators on $\mathfrak{h}$, $\vert\xi\rangle\langle
\chi\vert$ takes values in $ \mathfrak{L}_{1}( \mathfrak{h})$.
Applying Parseval's equality yields
%
%
\begin{equation}
\label{58}
\operatorname{tr} (A \vert\xi\rangle\langle\chi\vert)
=
\langle\chi, A \xi\rangle.
\end{equation}
Hence $\vert\xi\rangle\langle\chi\vert$ is $ \mathfrak{B}( \mathfrak
{L}_{1}( \mathfrak{h}) )$-measurable
because the dual of $ \mathfrak{L}_{1}( \mathfrak{h})$ is formed by all maps
$\varrho\mapsto\operatorname{tr} (A \varrho) $ with $A \in\mathfrak
{L}( \mathfrak{h})$.
Let $x,y \in\mathfrak{h}$.
The absolute value of the operator
$ \vert x \rangle\langle y \vert$
is equal to the operator
$
\vert y \rangle\langle y \vert
\Vert x \Vert/ \Vert y \Vert
$
in case $y \neq0$,
and coincides with the null operator otherwise.
Therefore
%
%
\begin{equation}
\label{55}
\Vert
\vert x \rangle\langle y \vert
\Vert_{1}
= \frac{\Vert x \Vert}{ \Vert y \Vert} \Vert y \Vert^2
= \Vert x \Vert\Vert y \Vert.
\end{equation}
Combining $\xi, \chi\in L^{2}( \mathbb{P},\mathfrak{h})$ with (\ref
{55}) gives
$\mathbb{E} \Vert\vert\xi\rangle\langle\chi\vert\Vert_{1} < \infty$,
and so
the Bochner integral $\mathbb{E}\vert\xi\rangle\langle\chi\vert$ is
well defined in the separable Banach space $\mathfrak{L}_{1}
( \mathfrak{h}) $.

We now turn to work in $\mathfrak{L}( \mathfrak{h}) $.
The application $( x, y ) \mapsto\vert x \rangle\langle y \vert$
from $\mathfrak{h} \times\mathfrak{h}$ to $\mathfrak{L}( \mathfrak
{h})$ is continuous, and in consequence
the measurability of $\xi$ and $\chi$ implies that
$\vert\xi\rangle\langle\chi\vert$ is $ \mathfrak{B} ( \mathfrak{L}(
\mathfrak{h}) )$-measurable.
Thus using
\mbox{$ \Vert\cdot\Vert_{\mathfrak{L}( \mathfrak{h})}
\leq\Vert\cdot\Vert_{1 }$}
we deduce that $\vert\xi\rangle\langle\chi\vert$ is Bochner
$\mathbb{P}$-integrable in $\mathfrak{L}( \mathfrak{h}) $; see, for
example,~\cite{Yosida1995} for a treatment of the Bochner integral in
Banach spaces which, in general, are not separable.
Since $ \mathfrak{L}_{1}( \mathfrak{h})$ is continuously embedded in
$\mathfrak{L}( \mathfrak{h}) $,
either of the interpretations of $\mathbb{E}\vert\xi\rangle\langle
\chi\vert$ given above refers to the same operator.

Finally,
for any $x,y$ belonging to $\mathfrak{h}$, the linear function $A
\mapsto\langle x, A y \rangle$ is continuous as a map from $\mathfrak
{L}( \mathfrak{h}) $ to $\mathbb{C}$.
This gives (\ref{53}).
Similarly, (\ref{58}) yields
$
\operatorname{tr} ( \mathbb{E}\vert\xi\rangle\langle\chi\vert)
=
\mathbb{E} \operatorname{tr} ( \vert\xi\rangle\langle\chi\vert)
= \mathbb{E} \langle\chi, \xi\rangle
$,
because $\operatorname{tr}( \cdot) \in\mathfrak{L}_{1}( \mathfrak
{h}) '$.
\end{pf}

\begin{remark}
Under the assumptions of Lemma~\ref{lema51}, $\mathbb{E}\vert\xi\rangle
\langle\chi\vert$ can also be interpreted as a Bochner integral in
the pointwise sense; see, for example,~\cite{DaPratoZabczyk1992}.
\end{remark}

To prove Theorem~\ref{teorema8},
we need the following lemma.

\begin{lemma}
\label{lema52}
Let $C$ be a self-adjoint positive operator in $\mathfrak{h}$.
Suppose that $\xi\in L^{2}_{C}( \mathbb{P},\mathfrak{h})$ and $A \in
\mathfrak{L}( ( \mathcal{D}( C) ,\Vert\cdot\Vert_{C}) ,\mathfrak{h}) $.
Then $A \xi$ belongs to $L^{2}( \mathbb{P},\mathfrak{h})$.
\end{lemma}

\begin{pf}
Since $A \xi= A \pi_{C} ( \xi)$ $\mathbb{P}$-a.s.,
from Remark~\ref{notaMedibilidad} we deduce that $A \xi$ is strongly
measurable.
Thus
$A \xi\in L^{2}( \mathbb{P},\mathfrak{h})$.
\end{pf}

\begin{pf*}{Proof of Theorem~\ref{teorema8}}
We start by proving statement (a).
Let $x \in\mathcal{D}( C)$ and let $y \in\mathfrak{h}$. Using Lemma
\ref{lema51} yields
\[
\langle Cx, \varrho y \rangle
=
\mathbb{E} \langle C x, \xi\rangle\langle\xi, y \rangle
=
\mathbb{E} \langle x, C \xi\rangle\langle\xi, y \rangle.
\]
In Lemma~\ref{lema52} we take $A=C$ to obtain $C \xi\in L^{2}( \mathbb
{P},\mathfrak{h})$. Thus, Lemma~\ref{lema51} implies
$
\mathbb{E} \langle x, C \xi\rangle\langle\xi, y \rangle
=
\langle x, \mathbb{E}\vert C \xi\rangle\langle\xi\vert y \rangle
$,
and so
$
\langle Cx, \varrho y \rangle
=
\langle x, \mathbb{E}\vert C \xi\rangle\langle\xi\vert y \rangle
$.
Then $\varrho y \in\mathcal{D}( C^{\ast}) = \mathcal{D}( C ) $ and
$
C\varrho y = \mathbb{E}\vert C \xi\rangle\langle\xi\vert y
$,
which is our assertion.

Part (a) yields
$
\mathcal{D}( B )
=
\mathcal{D}( A \varrho B ) $,
and so $A \varrho B$ is densely defined.
We next prove that
$A \varrho B$ coincides with
$ \mathbb{E}\vert A \xi\rangle\langle B^{\ast} \xi\vert$
on $\mathcal{D}( B )$.
For this purpose,
we approximate $A$ by $A R_{n}$, where $R_{n}$ is the Yosida
approximation of $-C$.

Suppose that $x \in\mathfrak{h}$ and $y \in\mathcal{D}( B )$. As in
the proof of Lemma~\ref{lema42} we consider
$
R_{n} = n (n+C )^{-1}
$,
and so
$C R_{n} z \longrightarrow_{n \rightarrow\infty} Cz$ for any $x \in
\mathcal{D}( C )$.
Therefore
$
\langle x, A R_{n} \varrho B y \rangle
\longrightarrow_{n \rightarrow\infty}
\langle x, A \varrho B y \rangle
$, and
hence
Lemma~\ref{lema51} gives
%
%
\begin{equation}
\label{56}\quad
\langle x, A \varrho B y \rangle
=
\lim_{n \rightarrow\infty} \mathbb{E} \langle( A R_{n} )^{\ast} x,
\xi\rangle\langle\xi, B y \rangle
=
\lim_{n \rightarrow\infty} \mathbb{E} \langle x, A R_{n} \xi\rangle
\langle\xi, B y \rangle.
\end{equation}
Since
$\Vert R_{n} \Vert\leq1$
and
$R_{n}$ commutes with $C$,
$\Vert A R_{n} z \Vert\leq K \Vert z \Vert_{C}$.
Using the dominated convergence theorem we obtain
%
%
\begin{equation}
\label{57}
\langle x, A \varrho B y \rangle
=
\lim_{n \rightarrow\infty} \mathbb{E} \langle x, A R_{n} \xi\rangle
\langle\xi, B y \rangle
=
\mathbb{E} \langle x, A \xi\rangle\langle B^{\ast}\xi, y \rangle.
\end{equation}

Since $B$ is densely defined,
$B^{\ast}$ is a closed operator.
Remark~\ref{nota1} now shows that
$B^{\ast} \in\mathfrak{L}( ( \mathcal{D}( C) ,\mbox{$\Vert\cdot\Vert_{C}$})
,\mathfrak{h})$,
and so
applying Lemma~\ref{lema52} gives
$
A \xi,\break B^{\ast} \xi\in L^{2}( \mathbb{P},\mathfrak{h})
$.
Combining (\ref{57}) with Lemma~\ref{lema51}
we get
$
\langle x, A \varrho B y \rangle
=\break
\langle x, \mathbb{E}\vert A \xi\rangle\langle B^{\ast} \xi\vert y
\rangle
$.
Since the closure of $A \varrho B$ is equal to $\mathbb{E}\vert A \xi
\rangle\langle B^{\ast} \xi\vert$,
we complete the proof of statement (b) by using Lemma~\ref{lema51}.
\end{pf*}

\subsection{\texorpdfstring{Proof of Theorem \protect\ref{teorema4}}{Proof of Theorem 3.1}}
First,
we easily construct a random variable that represents a given
$C$-regular operator.

\begin{lemma}
\label{lema6}
Let $\varrho\in\mathfrak{L}_{1,C}^{+}( \mathfrak{h}) $, with $C$
self-adjoint positive operator in $\mathfrak{h}$. Then there exists
$\xi\in L_{C}^{2}( \mathbb{P},\mathfrak{h}) $
such that $\varrho=\mathbb{E}\vert\xi\rangle\langle\xi\vert$
and $\Vert\xi\Vert^{2}=\operatorname{tr}( \varrho) $ a.s.
\end{lemma}

\begin{pf}
In case $\varrho=0$,
we take $\xi=0$.
Otherwise,
consider that $\varrho$ is written as in Definition~\ref{def2}.
Then, we choose $\Omega= \mathfrak{I}$,
and for any $n \in\mathfrak{I}$ we define
$\mathbb{P}( \{ n \} )= \lambda_{n}/\operatorname{tr}( \varrho)$
and $\xi( n ) = \sqrt{\operatorname{tr}( \varrho) }u_{n}$.
\end{pf}

Second,
we use part (a) of Theorem~\ref{teorema8},
together with Lemma~\ref{lema51},
to establish the sufficient condition of Theorem~\ref{teorema4}.

\begin{lemma}
\label{lema53}
Let $C$ be a self-adjoint positive operator in $\mathfrak{h}$. Suppose
that $\varrho= \mathbb{E}\vert\xi\rangle\langle\xi\vert$, with
$\xi\in L_{C}^{2}( \mathbb{P},\mathfrak{h}) $. Then $\varrho$ is $C$-regular.
\end{lemma}

\begin{pf}
Lemma~\ref{lema51} shows that $ \varrho\in\mathfrak{L}_{1}^{+}(
\mathfrak{h} )$,
hence
$
\varrho=\sum_{n\in\mathfrak{I}}\lambda_{n}\vert u_{n}\rangle\langle
u_{n}\vert
$,
where $\mathfrak{I}$ is a countable set, $( \lambda_{n}) _{n\in\mathfrak
{I}}$ are summable positive real numbers and $( u_{n}) _{n\in\mathfrak
{I}}$ is a orthonormal family of vectors of $\mathfrak{h}$.
Using statement (a) of Theorem~\ref{teorema8}
yields $ u_{n} \in\mathcal{D}( C) $ for all $n\in\mathfrak{I}$.

We can extend $( u_{n}) _{n\in\mathfrak{I}}$ to an orthonormal basis
$( e_{n}) _{n\in\mathfrak{I}'}$ of $\mathfrak{h}$ formed by elements of
$\mathcal{D}( C) $.
From Parseval's equality we obtain
\[
\sum_{n\in\mathfrak{I}}\lambda_{n}\Vert Cu_{n}\Vert^{2}
=
\sum_{n\in\mathfrak{I}}\sum_{k\in\mathfrak{I}'}\lambda_{n}\vert
\langle Cu_{n},e_{k}\rangle\vert^{2}
=
\sum_{k\in\mathfrak{I}'}\sum_{n\in\mathfrak{I}}\lambda_{n}\langle
Ce_{k},\vert u_{n}\rangle\langle u_{n}\vert Ce_{k}\rangle,
\]
and so
$
\sum_{n\in\mathfrak{I}}\lambda_{n}\Vert Cu_{n}\Vert^{2}
=
\sum_{k\in\mathfrak{I}'}\langle Ce_{k},\varrho Ce_{k}\rangle
$.
Combining Lemma~\ref{lema51} with Parseval's equality we now get
\[
\sum_{n\in\mathfrak{I}}\lambda_{n}\Vert Cu_{n}\Vert^{2}
=
\sum_{k\in\mathfrak{I}'}\mathbb{E}\vert\langle
\xi,Ce_{k}\rangle\vert^{2}
=
\mathbb{E}\sum_{k\in\mathfrak{I}'}\vert\langle C\xi
,e_{n}\rangle\vert^{2}
=
\mathbb{E}\Vert C\xi\Vert^{2}.
\]
This gives $\varrho\in\mathfrak{L}_{1,C}^{+}( \mathfrak{h}) $.
\end{pf}

\subsection{\texorpdfstring{Proof of Theorem \protect\ref{teor7}}{Proof of Theorem 4.1}}
\label{subsecteor7}

We first establish, in our framework, the well-known relation between
Heisenberg and Schr\"{o}dinger pictures.

\begin{lemma}
\label{lema10}
Suppose that Hypothesis~\ref{HipN5} holds,
together with $\xi\in L_{C}^{2}( \mathbb{P},\mathfrak{h}) $.
Let $\mathcal{T}_{t}( A) $ be as in Theorem~\ref{teorema3}.
Then for all $A\in\mathfrak{L}( \mathfrak{h}) $,
%
%
\begin{equation}
\label{35}
\operatorname{tr}( A\mathbb{E}\vert X_{t}( \xi) \rangle\langle X_{t}(
\xi) \vert) =\operatorname{tr}( \mathcal{T}_{t}(
A) \mathbb{E}\vert\xi\rangle\langle\xi\vert).
\end{equation}
\end{lemma}

\begin{pf}
Fix $A\in\mathfrak{L}( \mathfrak{h}) $,
and
define the function $f_{n} \dvtx \mathfrak{h} \rightarrow\mathbb{C}$ by
$
f_{n}( x) = \langle x,Ax\rangle
$
if $\Vert x\Vert\leq n$,
and
$
f_{n}( x) = 0
$
otherwise.
Using the Markov property of $X_t ( \xi)$,
which can be obtained by techniques of well-posed martingale problems,
we get
%
%
\begin{equation}
\label{n311}
\mathbb{E}( f_{n}( X_{t}( \xi) ) )
=
\mathbb{E}
\bigl(
( f_{n}( X_{t}( \xi) ) )
\diagup\mathfrak{F}_{0}
\bigr)
=
\mathbb{ E}P_{t}f_{n}( \xi),
\end{equation}
where
$P_{t}f_{n}( x) =\mathbb{E}( f_{n}( X_{t}( x) ) ) $
for all $x\in\mathcal{D}( C) $.

We will take the limit as $n \rightarrow\infty$ in (\ref{n311}).
The dominated convergence theorem leads to
%
%
\begin{equation}
\label{62}
\lim_{n \rightarrow\infty} \mathbb{E}( f_{n}( X_{t}( \xi) ) )
=
\mathbb{E}\langle X_{t}( \xi) ,AX_{t}( \xi) \rangle.
\end{equation}
Combining (\ref{62}) with (\ref{42}) yields
$
P_{t}f_{n}( x) \longrightarrow_{n\rightarrow\infty}
\langle x,\mathcal{T}_{t}( A) x\rangle
$
whenever $x\in\mathcal{D}( C) $.
Since $\Vert P_{t}f_{n}( x) \Vert\leq\Vert A \Vert\Vert x \Vert^{2}$,
according to the dominated convergence theorem, we have
$
\mathbb{ E}P_{t}f_{n}( \xi)
\longrightarrow
\mathbb{ E} \langle\xi,\mathcal{T}_{t}( A) \xi\rangle
$
as
$
n \rightarrow\infty
$.
Then,
letting $n \rightarrow\infty$ in (\ref{n311})
we get
$
\mathbb{E}\langle X_{t}( \xi) ,AX_{t}( \xi)
\rangle=\mathbb{E}\langle\xi,\mathcal{T}_{t}( A)
\xi\rangle
$
by (\ref{62}),
and so
Theorem~\ref{teorema8} leads to (\ref{35}).
\end{pf}

We next check that
$\rho_{t}( \varrho) $ is well defined by (\ref{31}).

\begin{lemma}
\label{lema13}
Let Hypothesis~\ref{HipN5} hold
and consider
$\xi, \varphi\in L_{C}^{2}( \mathbb{P},\mathfrak{h}) $
such that
$\mathbb{E}\vert\xi\rangle\langle\xi\vert=\mathbb{E}\vert
\varphi\rangle\langle\varphi\vert$.
Then
$
\mathbb{E}\vert X_{t}( \xi) \rangle\langle X_{t}(
\xi) \vert=\mathbb{E}\vert X_{t}( \varphi)
\rangle\langle X_{t}( \varphi) \vert
$.
\end{lemma}

\begin{pf}
Let $A\in\mathfrak{L}( \mathfrak{h}) $. Using Lemma~\ref{lema10}
yields
\[
\operatorname{tr}( A\mathbb{E}\vert X_{t}( \xi) \rangle\langle X_{t}(
\xi) \vert)
=
\operatorname{tr}( \mathcal{T}_{t}( A) \mathbb{E}\vert\xi\rangle
\langle\xi\vert)
=
\operatorname{tr}( A\mathbb{E}\vert X_{t}( \varphi) \rangle\langle
X_{t}( \varphi) \vert) .
\]
Hence $\Vert\mathbb{E}\vert X_{t}( \xi) \rangle\langle
X_{t}( \xi) \vert-\mathbb{E}\vert X_{t}(
\varphi) \rangle\langle X_{t}( \varphi) \vert
\Vert_{\mathfrak{L}_{1}( \mathfrak{h}) }=0$; see, for example,
Proposition 9.12 of~\cite{Parthasarathy1992}.
\end{pf}

We now address the contraction property of the restriction of $\rho_{t}$
to $\mathfrak{L}_{1,C}^{+} ( \mathfrak{h}) $.

\begin{lemma}
\label{lema11}
Let Hypothesis~\ref{HipN5} hold. If $\varrho,\widetilde{\varrho}
$ are $C$-regular, then
%
%
\begin{equation}
\label{n312}
{\operatorname{tr}}\vert\rho_{t}( \varrho) -\rho_{t}( \widetilde{%
\varrho}) \vert\leq{\operatorname{tr}}\vert\varrho-\widetilde{\varrho}
\vert.\vadjust{\goodbreak}
\end{equation}
\end{lemma}

\begin{pf}
Since
$
{\operatorname{tr} }\vert\rho_{t}( \varrho) -
\rho_{t}( \widetilde{\varrho}) \vert
=
\sup_{\Vert A\Vert_{\mathfrak{L}( \mathfrak{h} )} =1}
\vert{\operatorname{tr}}( A\rho_{t}(
\varrho) ) -\operatorname{tr}( A\rho_{t}( \widetilde{\varrho})
) \vert
$,
according to Lemma~\ref{lema10} we have
\[
{\operatorname{tr}}\vert\rho_{t}( \varrho) -\rho_{t}( \widetilde{\varrho
}) \vert
=
\sup_{A\in\mathfrak{L}( \mathfrak{h}) ,\Vert A\Vert=1}
\vert{\operatorname{tr}}( \mathcal{T}_{t}( A) \varrho)
-\operatorname{tr}( \mathcal{T}_{t}( A) \widetilde{\varrho})
\vert.
\]
Therefore
$
{\operatorname{tr}}\vert\rho_{t}( \varrho) -\rho_{t}( \widetilde{\varrho
}) \vert
\leq
{\operatorname{tr}}\vert\varrho-\widetilde{\varrho}\vert
\sup_{\Vert A\Vert_{\mathfrak{L}( \mathfrak{h})} =1}
\Vert\mathcal{T}_{t}( A) \Vert
$,
and so Theorem~\ref{teorema3} leads to (\ref{n312}).
\end{pf}

The following lemma helps us to extend $\rho_{t}$ to all $\mathfrak
{L}_{1}( \mathfrak{h}) $.

\begin{lemma}
\label{lema12}
Suppose that $C$ is a self-adjoint positive operator in $\mathfrak{h}$.
Then $\mathfrak{L}_{1,C}^{+}( \mathfrak{h}) $ is dense in
$\mathfrak{L}_{1}^{+}( \mathfrak{h}) $ with respect to the trace norm.
\end{lemma}

\begin{pf}
Let $\varrho\in\mathfrak{L}_{1}^{+}( \mathfrak{h}) $.
Then there exists a sequence of orthonormal vectors $( u_{j}) _{j\in
\mathbb{N}}$
for which $\varrho=\sum_{j\in\mathbb{N}}\lambda_{j}\vert u_{j}\rangle
\langle u_{j}\vert$,
with $\lambda_{j}\geq0$ and $\sum_{j\in\mathbb{N}}\lambda_{j}<\infty$.
For any $x,y\in\mathfrak{h}$ we have
\[
{\operatorname{tr}}\bigl\vert
\vert x\rangle\langle x\vert-\vert
y\rangle\langle y\vert
\bigr\vert
=
\sup_{\Vert A\Vert_{\mathfrak{L}( \mathfrak{h})} =1}
\vert\langle
x,Ax\rangle-\langle y,Ay\rangle\vert
\leq
\Vert x-y\Vert^{2}+2\Vert y\Vert\Vert x-y\Vert,
\]
and so
$\{ \vert x\rangle\langle x\vert\dvtx x\in\mathcal{D}(
C) \} $ is a \mbox{$\Vert\cdot\Vert_{\mathfrak{L}%
_{1}( \mathfrak{h}) }$}-dense subset of $\{ \vert
x\rangle\langle x\vert\dvtx x\in\mathfrak{h}\} $
since
$\mathcal{D}( C) $ is dense in $\mathfrak{h}$.
Now,
the lemma follows from
$
{\operatorname{tr}}\vert\varrho-\sum_{j=1}^{n}\lambda_{j}\vert
u_{j}\rangle\langle u_{j}\vert\vert
=
\sum_{j=n+1}^{\infty}\lambda_{j}
\longrightarrow_{n\rightarrow\infty} 0
$.
\end{pf}

\begin{pf*}{Proof of Theorem~\ref{teor7}}
Combining Theorem~\ref{teorema4} with Lemma~\ref{lema13}
we obtain that (\ref{31}) defines unambiguously a linear
operator $\rho_{t}( \varrho) $
for any $\varrho\in\mathfrak{L}_{1,C}^{+}( \mathfrak{h}) $ and $t\geq0$.
Lemma~\ref{lema12} guarantees the uniqueness of the operator
belonging to $\mathfrak{L}( \mathfrak{L}_{1}( \mathfrak{h})
) $ for which (\ref{31}) holds.
We next extend $\rho_{t}$ to a bounded linear operator in $\mathfrak
{L}_{1}( \mathfrak{h}) $
by means of density arguments.

Suppose that $\varrho\in\mathfrak{L}_{1}^{+}( \mathfrak{h}) $.
By Lemma~\ref{lema12}, there exists a sequence $(
\varrho_{n}) _{n\in\mathbb{N}}$ of $C$-regular operators for which
$
\lim_{n\rightarrow\infty}
\Vert\varrho-\varrho_{n}\Vert_{\mathfrak{L}_{1}( \mathfrak{%
h}) }\rightarrow0
$.
We define
$
\rho_{t}( \varrho)
$
to be the limit in
$\mathfrak{L}_{1}( \mathfrak{h} ) $ of $\rho_{t}( \varrho_{n}) $ as
$n\rightarrow\infty$;
according to Lemma~\ref{lema11}
this limit exists and
does not depend on the choice of $( \varrho_{n}) _{n\in\mathbb{N}}$.
Recall that every $A \in\mathfrak{L}( \mathfrak{h}) $
has a unique decomposition of the form $A=\Re( A) +i$ $%
\Im( A) $,
with $\Re( A) $ and $\Im( A) $ self-adjoint operators in $\mathfrak{h}$.
For each $\varrho\in\mathfrak{L}_{1}( \mathfrak{h}) $ we set
\[
\rho_{t}( \varrho) =\rho_{t}( \Re( \varrho)
_{+}) -\rho_{t}( \Re( \varrho) _{-})
+i\bigl( \rho_{t}( \Im( \varrho) _{+})
-\rho_{t}( \Im( \varrho) _{-}) \bigr) ,
\]
where $A_{+}$, $A_{-}$ denotes, respectively, the positive and negative parts
of the self-adjoint operator $A$; see, for example, \cite
{BratteliRobinson1987} for
details.

We will verify that
$\rho_{t} \in\mathfrak{L} ( \mathfrak{L}_{1}( \mathfrak{h}) ) $.
Let $\varrho=\varrho_{1}-\varrho_{2}+i( \varrho_{3}-\varrho_{4}) $,
with $\varrho_{j}\in\mathfrak{L}_{1,C}^{+}( \mathfrak{h} ) $ for any
$j=1,\ldots,4$.
Since
$\Vert\mathcal{T}_{t}( A) \Vert\leq\Vert A \Vert$,
Lemma~\ref{lema10} yields
\[
{\operatorname{tr}}\vert\rho_{t}( \varrho) \vert
=
{\sup_{\Vert A\Vert_{\mathfrak{L}( \mathfrak{h}) } =1}}
\vert{\operatorname{tr}}( A\rho_{t}( \varrho) ) \vert
=
{\sup_{\Vert A\Vert_{\mathfrak{L}( \mathfrak{h}) } =1}}
\vert{\operatorname{tr}}( \mathcal{T}_{t}( A) \varrho) \vert
\leq
\operatorname{tr}( \vert\varrho\vert) .
\]
The construction of $\rho_{t}( \varrho) $ now implies
$
\Vert\rho_{t}( \varrho) \Vert_{\mathfrak{L}_{1}( \mathfrak{h}) }\leq
\Vert\varrho\Vert_{\mathfrak{L}_{1}( \mathfrak{h}) }
$
for all $\varrho\in\mathfrak{L}_{1}( \mathfrak{h}) $.
Consider two $C$-regular operators $\varrho,\widetilde{\varrho}$ and
$\alpha\geq0$.
By Definition~\ref{def2},\vadjust{\goodbreak} $\varrho+\alpha\widetilde{\varrho}$ belongs
to $\mathfrak{L}_{1,C}^{+}( \mathfrak{h}) $.
If $A\in\mathfrak{L}( \mathfrak{h}) $,
then applying Lemma~\ref{lema10} we obtain
\[
\operatorname{tr}\bigl( \rho_{t}( \varrho+\alpha\widetilde{\varrho}) A\bigr)
=
\operatorname{tr}( \mathcal{T}_{t}( A) \varrho)
+
\alpha\operatorname{tr}( \mathcal{T}_{t}( A) \widetilde{\varrho}) \\
=
\operatorname{tr}\bigl( \bigl( \rho_{t}( \varrho)
+
\alpha
\rho_{t}( \widetilde{\varrho}) \bigr) A\bigr) .
\]
Therefore $\Vert\rho_{t}( \varrho+\alpha\widetilde{\varrho}%
) -\rho_{t}( \varrho) -\alpha\rho_{t}( \widetilde{%
\varrho}) \Vert_{\mathfrak{L}_{1}( \mathfrak{h}) }=0$%
, and so
Lemma~\ref{lema12} leads to $\rho_{t}( \varrho+\alpha%
\widetilde{\varrho}) =\rho_{t}( \varrho)
+\alpha\rho_{t}( \widetilde{\varrho}) $ for any $\varrho,%
\widetilde{\varrho}\in\mathfrak{L}_{1}^{+}( \mathfrak{h}) $.
Careful algebraic manipulations now show the linearity of $\rho_{t}\dvtx%
\mathfrak{L}_{1}( \mathfrak{h}) \rightarrow\mathfrak{L}%
_{1}( \mathfrak{h}) $.
\end{pf*}

\subsection{\texorpdfstring{Proof of Theorem \protect\ref{teor9}}{Proof of Theorem 4.2}}
\label{subsecteor9}

Let us first prove the continuity of the map
$
\xi\mapsto\rho_{t}( \mathbb{E}\vert
\xi\rangle\langle\xi\vert)
$.

\begin{lemma}
\label{lema16}
Assume that Hypothesis~\ref{HipN5} holds.
Let $\xi$ and $\xi_{n}$, with \mbox{$n\in\mathbb{N}$}, be random variables in
$L_{C}^{2}( \mathbb{P},\mathfrak{h}) $
satisfying
$\mathbb{E}\Vert\xi-\xi_{n}\Vert^{2}\longrightarrow_{n\rightarrow
\infty}0$.
Then\break
$\rho_{t}( \mathbb{E}\vert\xi_{n}\rangle\langle\xi_{n}\vert) $
converges in $\mathfrak{L}( \mathfrak{h} ) $ to
$\rho_{t}( \mathbb{E}\vert\xi\rangle\langle\xi\vert) $
as $n \rightarrow\infty$.
\end{lemma}

\begin{pf}
Let $x \in\mathfrak{h}$. Combining (\ref{31}) with the linearity of
(\ref{2}) we get
\begin{eqnarray*}
&&
\bigl\Vert\rho_{t}( \mathbb{E}\vert\xi_{n}\rangle\langle\xi_{n}\vert)
x-\rho_{t}( \mathbb{E}\vert\xi
\rangle\langle\xi\vert) x\bigr\Vert
\\
&&\qquad \leq
\mathbb{E}\vert\langle X_{t}( \xi_{n}) ,x\rangle\vert\Vert X_{t}( \xi
_{n}- \xi) \Vert
+
\mathbb{E}\vert\langle X_{t}( \xi- \xi_{n}) ,x\rangle\vert\Vert
X_{t}( \xi) \Vert
\\
&&\qquad \leq
\Vert x\Vert
\bigl(
\mathbb{E}\Vert\xi-\xi_{n}\Vert^{2} +
2 \sqrt{\mathbb{E}\Vert\xi-\xi_{n}\Vert^{2}} \sqrt{ \mathbb{E}\Vert\xi
\Vert^{2}}
\bigr).
\end{eqnarray*}
In the last inequality we used that
$
\mathbb{E}\Vert X_{t}( \eta) \Vert^{2}
\leq
\mathbb{E} \Vert\eta\Vert^{2}
$
for $\eta\in L_{C}^{2}( \mathbb{P},\mathfrak{h}) $.
\end{pf}

\begin{pf*}{Proof of Theorem~\ref{teor9}}
There exits a sequence $( \xi_{n}) _{n}$ of
\mbox{$( \mathcal{D}( C), \Vert\cdot\Vert)$}-valued random variables with
finite ranges
such that $\Vert\xi_{n}-\xi\Vert$ converges monotonically to $0$;
see, for example,~\cite{DaPratoZabczyk1992}.
By Lemma~\ref{lema16},
$
\rho_{t}( \mathbb{E}\vert\xi_{n}\rangle\langle\xi_{n}\vert)
$
converges to
$
\rho_{t}( \mathbb{E}\vert\xi\rangle\langle\xi\vert)
$
in $\mathfrak{L}( \mathfrak{h})$.
Since $\rho_{t}$ is linear, an easy computation shows that
$
\mathbb{E}\rho_{t}( \vert\xi_{n}\rangle\langle\xi_{n}\vert
)
=
\rho_{t}( \mathbb{E}\vert\xi_{n}\rangle
\langle\xi_{n}\vert)
$,
hence
%
%
\begin{equation}
\label{n37}
\mathbb{E}\rho_{t}( \vert\xi_{n}\rangle\langle\xi_{n}\vert
)
\longrightarrow_{n\rightarrow\infty}
\rho_{t}( \mathbb{E}\vert\xi\rangle\langle\xi\vert)
\qquad\mbox{in } \mathfrak{L}( \mathfrak{h}) .
\end{equation}

We will prove that
$\mathbb{E}\rho_{t}( \vert\xi_{n}\rangle\langle\xi_{n}\vert) $
converges to
$\mathbb{E}\rho_{t}( \vert\xi\rangle\langle\xi\vert) $ in $\mathfrak
{L}( \mathfrak{h}
) $
as $n \rightarrow\infty$,
which together with (\ref{n37})
implies
$
\rho_{t}( \mathbb{E}\vert\xi\rangle\langle\xi\vert)
=
\mathbb{E} \rho_{t}( \vert\xi\rangle\langle\xi\vert)
$.
From Lem\-ma~\ref{lema16} we obtain
$
\Vert\rho_{t}( \vert\xi_{n}\rangle\langle\xi_{n}\vert
) -\rho_{t}( \vert\xi\rangle\langle\xi\vert)
\Vert_{\mathfrak{L}( \mathfrak{h}) }\longrightarrow
_{n\rightarrow\infty}0
$.
For any $x,y \in\mathfrak{h}$ we have
$
\Vert
\vert x \rangle\langle y \vert
\Vert_{1}
= \Vert x \Vert\Vert y \Vert
$,
and so Lemma~\ref{lema11} yields
\[
\Vert\rho_{t}( \vert\xi_{n}\rangle\langle\xi_{n}\vert
) \Vert
\leq
\Vert\rho_{t}( \vert\xi_{n}\rangle\langle\xi_{n}\vert
) \Vert_{1}
\leq
\Vert\xi_{n}\Vert^{2}
\leq
2( \Vert\xi_{1}-\xi\Vert^{2}+\Vert
\xi\Vert^{2}) .
\]
Therefore
$
\mathbb{E}\Vert\rho_{t}( \vert\xi_{n}\rangle\langle\xi
_{n}\vert) -\rho_{t}( \vert\xi\rangle\langle
\xi\vert) \Vert_{\mathfrak{L}( \mathfrak{h})
}\longrightarrow_{n\rightarrow\infty}0.
$
\end{pf*}

\subsection{\texorpdfstring{Proof of Theorem \protect\ref{teor8}}{Proof of Theorem 4.3}}
Our proof is divided into three lemmata.
The first two deal with the semigroup property of $( \rho_t )_{t \geq
0} $.

\begin{lemma}
\label{lema14}
Let Hypothesis~\ref{HipN5} hold,
and
let $\varrho$ be $C$-regular.
Then
for all $t \geq0$,
$\rho_{t}( \varrho)$ belongs to $\mathfrak{L}_{1,C}^{+}( \mathfrak{h}) $
and
$\rho_{t+s}( \varrho) =\rho_{t}\circ\rho_{s}( \varrho) $
whenever $s\geq0$.
\end{lemma}

\begin{pf}
Since
$X_{t}( \xi) \in L_{C}^{2}( \mathbb{P},\mathfrak{h}) $,
combining Theorem~\ref{teorema4} with (\ref{31}) gives
$
\rho_{t}( \mathfrak{L}_{1,C}^{+}( \mathfrak{h}) )
\subset\mathfrak{L}_{1,C}^{+}( \mathfrak{h})
$.\vadjust{\goodbreak}

We will establish the semigroup property of the restriction of $ \rho$ to
$ \mathfrak{L}_{1,C}^{+}( \mathfrak{h})$.
Consider
$\xi\in L_{C}^{2}( \mathbb{P},\mathfrak{h}) $ satisfying
$\rho=\mathbb{E}\vert\xi\rangle\langle\xi\vert$,
and fix
$ x,y \in\mathfrak{h}$.
For all $z \in\mathfrak{h}$
we define
$p_{n}( z )
=
\langle z,x\rangle\langle y,z\rangle
$
if
$
\vert\langle z,x\rangle\langle y,z\rangle\vert\leq n
$,
and
$p_{n}( z ) = 0$
otherwise.
Using the Markov property of
$
X_{t}( \xi)
$
we deduce that
%
%
\begin{equation}
\label{38}
\mathbb{E}( p_{n} ( X_{t+s}( \xi) ) )
=
\mathbb{E} \bigl(
( p_{n} ( X_{t+s}( \xi) ) )
\diagup\mathfrak{F}_{s}
\bigr)
=
\mathbb{E}P_{t}( p_{n} ) ( X_{s}( \xi)
) ,
\end{equation}
where for all $z \in\mathcal{D}( C) $,
$P_{t}( p_{n} ) ( z) = \mathbb{E}( p_{n} ( X_{t}( z) ) ) $.

Let $z \in\mathcal{D}( C) $.
Applying the dominated convergence theorem gives
\[
\lim_{n\rightarrow\infty} \mathbb{E}( p_{n} ( X_{t}( z) ) )
=
\mathbb{E} \langle X_{t}( z) ,x\rangle\langle y,X_{t}( z) \rangle
=
\langle y,\rho_{t}( \vert z\rangle\langle z\vert) x\rangle,
\]
hence
$
\lim_{n\rightarrow\infty} P_{t}( p_{n} ) ( z)
=
\langle y,\rho_{t}( \vert z\rangle\langle z\vert) x\rangle
$.
Then
$
\mathbb{E} P_{t} ( p_{n} ) ( X_{s} ( \xi) )
\longrightarrow_{n\rightarrow\infty}\break
\mathbb{E}
\langle y,\rho_{t}( \vert X_{s}( \xi)
\rangle\langle X_{s}( \xi) \vert) x\rangle
$,
and so
Theorem~\ref{teor9} leads to
%
%
\begin{equation}
\label{331}\quad\qquad
\lim_{n\rightarrow\infty}
\mathbb{E} P_{t} ( p_{n} ) ( X_{s} ( \xi) )
=
\langle y,\rho_{t}( \mathbb{E}\vert X_{s}( \xi)
\rangle\langle X_{s}( \xi) \vert) x\rangle
=
\langle y,\rho_{t} \circ\rho_{s} ( \varrho) x\rangle.
\end{equation}
By (\ref{38}),
in (\ref{331}) we replace
$s$ by $0$
and
$t$ by $t+s$
to obtain
\[
\lim_{n\rightarrow\infty}
\mathbb{E}( p_{n} ( X_{t+s}( \xi) ) )
=
\lim_{n\rightarrow\infty}
\mathbb{E} P_{t+s} ( p_{n} ) ( X_{0} ( \xi) )
=
\langle y,\rho_{t+s}( \varrho) x\rangle.
\]
Thus,
letting $n\rightarrow\infty$ in (\ref{38}),
we get
$\rho_{t+s}( \varrho) =\rho_{t}\circ\rho_{s}( \varrho) $
by (\ref{331}).
\end{pf}

\begin{lemma}
\label{lema24}
Under Hypothesis~\ref{HipN5}, $( \rho_{t})_{t\geq0}$ is a semigroup of
contractions
which leaves $\mathfrak{L}_{1}^{+}( \mathfrak{h}) $ invariant.
\end{lemma}

\begin{pf}
By Theorem~\ref{teor7},
$\Vert\rho_{t}\Vert_{\mathfrak{L}( \mathfrak{L}_{1}( \mathfrak{h}) )
}\leq1$.
Since $\rho_{t}( \varrho) $ is positive whenever $\varrho$ is $C$-regular,
using Lemma~\ref{lema12} yields
$\langle x,\rho_{t}( \varrho) x\rangle\geq0$
for any $\varrho\in\mathfrak{L}_{1}^{+}( \mathfrak{h}) $ and $x\in
\mathfrak{h}$.

Suppose that
$\varrho=\varrho_{1}-\varrho_{2}+i( \varrho_{3}-\varrho_{4}) $,
where $\varrho_{1},\ldots,\varrho_{4}$ are $C$-regular operators.
Applying (\ref{31}) gives $\rho_{0}( \varrho) =\varrho$,
and Lemma~\ref{lema14} asserts that
$\rho_{t+s}( \varrho) =\rho_{t}\circ\rho_{s}( \varrho) $ for any
$s,t\geq0$.
Then,
combining Lemma~\ref{lema12} with density arguments, we deduce that
$( \rho_{t}) _{t\geq0}$ is a semigroup.
\end{pf}

We now examine the continuity of the map
$t\mapsto$ $\rho_{t}( \varrho) $ when $\varrho$ is $C$-regular.

\begin{lemma}
\label{lema25}
Adopt Hypothesis~\ref{HipN5}, together with $\varrho\in%
\mathfrak{L}_{1,C}^{+}( \mathfrak{h}) $. Then the map $%
t\mapsto\rho_{t}( \varrho) $ from $[ 0,\infty[ $ to $%
\mathfrak{L}_{1}( \mathfrak{h}) $ is continuous.
\end{lemma}

\begin{pf}
Consider $\xi\in L_{C}^{2}( \mathbb{P},\mathfrak{h}) $ such that $%
\varrho=\mathbb{E}\vert\xi\rangle\langle\xi\vert$.
Theorem~\ref{teorema8} yields
$
\mathbb{E}\Vert X_{t}( \xi) \Vert^{2}
\leq
\mathbb{E}\Vert\xi\Vert^{2}
=
\operatorname{tr}( \varrho)
$
for all $t\geq0$,
and so combining Theorem~\ref{teorema8} with the Cauchy--Schwarz
inequality yields
\begin{eqnarray*}
{\operatorname{tr}}\vert\rho_{t}( \varrho) -\rho_{s}( \varrho)
\vert &=& \sup_{A\in\mathfrak{L}( \mathfrak{h}) ,\Vert
A\Vert=1}\vert\mathbb{E}\langle X_{t}( \xi)
,AX_{t}( \xi) \rangle-\langle X_{s}( \xi)
,AX_{s}( \xi) \rangle\vert\\
&\leq& 2( \operatorname{tr}( \varrho) ) ^{1/2}\bigl( \mathbb{E}%
\Vert X_{t}( \xi) -X_{s}( \xi) \Vert
^{2}\bigr) ^{1/2}.
\end{eqnarray*}
Since
$
\mathbb{E} ( \sup_{s \in[0, T ]} \Vert X_s( \xi) \Vert^{2} )
< \infty
$
for any $T>0$
(see, e.g., Theorem 4.2.5 of~\cite{Prevot2007}),
using the dominated convergence theorem,
we get (\ref{310}).
\end{pf}

\subsection{\texorpdfstring{Proof of Theorem \protect\ref{teor10}}{Proof of Theorem 4.4}}
\label{subsecteor10}

First, we establish the weak continuity of the map $t\mapsto AX_{t}( \xi
) $ when $A$ is relatively bounded by $C$.

\begin{lemma}
\label{lema17}
Assume that Hypothesis~\ref{HipN5} holds.
If $\xi$ belongs to $L_{C}^{2}( \mathbb{P},\mathfrak{h}) $
and if $A \in\mathfrak{L}( ( \mathcal{D}( C) ,\mbox{$\Vert\cdot\Vert_{C}$})
,\mathfrak{h}) $,
then
for all $\psi\in L^{2}( \mathbb{P},\mathfrak{h}) $ and $t\geq0$ we have
%
%
\begin{equation} \label{n39}
\lim_{s\rightarrow t}\mathbb{E}\langle\psi,AX_{s}( \xi)
\rangle=\mathbb{E}\langle\psi,AX_{t}( \xi)
\rangle.
\end{equation}
\end{lemma}

\begin{pf}
Let $( s_{n}) _{n}$ be a sequence of nonnegative
real numbers converging to~$t$.
Since
$( ( X_{s_{n}}( \xi) ,AX_{s_{n}}( \xi) ,CX_{s_{n}}( \xi) ) )_{n}$
is a bounded sequence in $L^{2}( \mathbb{P}, \mathfrak{h}^{3} ) $
with $\mathfrak{h}^{3} = \mathfrak{h}\times\mathfrak{h}\times
\mathfrak{h}$,
there exists a subsequence
$( s_{n( k) }) _{k}$ for which
%
%
\begin{equation}
\label{n38} \bigl( X_{s_{n( k) }}( \xi) ,AX_{s_{n( k)}}( \xi)
,CX_{s_{n( k) }}( \xi) \bigr) \longrightarrow_{k\rightarrow\infty} (
Y,U,V)
\end{equation}
weakly in $L^{2}( \mathbb{P},\mathfrak{h}^{3})$.

Set
$\mathfrak{M}=\{ ( \eta,A\eta,C\eta) \dvtx\eta\in L_{C}^{2}( \mathbb
{P},\mathfrak{h}) \} $.
Then $\mathfrak{M}$ is a linear manifold of $L^{2}( \mathbb{P},\mathfrak
{h}^{3}) $ closed with respect to the strong topology.
In fact, suppose that\vspace*{1pt}
$( ( \eta_{n},A\eta_{n}, C\eta_{n}) ) _{n}$
is a sequence of elements of $\mathfrak{M}$ that converges to $( \eta
_{1},\eta_{2},\eta_{3}) $ in $L^{2}(
\mathbb{P},\mathfrak{h}^{3}) $.
Hence there exists a subsequence
$( ( \eta_{n( j) },A\eta_{n( j)}, C\eta_{n( j) }) ) _{j}$
converging almost surely to $( \eta_{1},\eta_{2},\eta_{3}) $.
Therefore $\eta_{1} \in\mathcal{D}( C )$ and $\eta_{3}=C\eta_{1}$ by
$C$ is closed.
Using $A\in\mathfrak{L}( ( \mathcal{D}( C) ,\mbox{$\Vert\cdot\Vert_{C}$})
,\mathfrak{h}) $ gives $\eta_{2}=A\eta_{1}$.

For any $k \in\mathbb{N}$,
$(X_{s_{n( k) }}( \xi) ,AX_{s_{n( k) }}( \xi) ,CX_{s_{n( k) }}( \xi) )
$ belongs to $\mathfrak{M}$.
Since
$\mathfrak{M}$ is a closed linear manifold of $L^{2}( \mathbb
{P},\mathfrak{h}^{3}) $,
(\ref{n38}) implies
$( Y,U,V) \in\mathfrak{M}$; see, for example, Section III.1.6 of
\cite{Kato1995}.
Combining
the dominated convergence theorem
with
$
\mathbb{E} ( \sup_{s \in[0, t+1 ]} \Vert X_s( \xi) \Vert^{2} )
< \infty
$
we get that
$
\mathbb{E} \Vert X_{s_{n( k)}}( \xi) -X_{t}( \xi) \Vert^{2}
$
converges\vspace*{1pt} to $0$.
Thus
$Y=X_{t}( \xi) $, and so
$U=AX_{t}( \xi) $.
Hence
$
AX_{s_{n( k) }}( \xi)
$
converges to
$AX_{t} ( \xi) $
weakly in $L^{2}( \mathbb{P},\mathfrak{h})$.
\end{pf}

Second,
we show that the probabilistic representation of
the right-hand side of (\ref{312}) is continuous as a function from $[
0,+\infty[ $ to $\mathbb{C}$.

\begin{lemma}
\label{lema18}
Let Hypothesis~\ref{HipN5} hold.
Fix $\xi\in L_{C}^{2}( \mathbb{P},\mathfrak{h}) $ and $A \in\mathfrak
{L}( \mathfrak{h}) $.
Then,
the function that maps each $t$ in $[ 0,+\infty[ $
to the complex number
$
\mathbb{E}\langle GX_{t}( \xi)$, $AX_{t}( \xi) \rangle
+
\mathbb{E}\langle X_{t}( \xi) ,AG X_{t}( \xi) \rangle
+\sum_{k=1}^{\infty}
\mathbb{E}\langle L_{k}X_{t}( \xi) ,AL_{k}X_{t}( \xi) \rangle
$
is continuous.
\end{lemma}

\begin{pf}
Let $( t_{n}) _{n}$ be a sequence of nonnegative real numbers such that
$t_{n}$ converges to $t$.
Since
$
\mathbb{E} ( \sup_{s \in[0, t+1 ]} \Vert X_s( \xi) \Vert^{2} )
< \infty
$
(see, e.g., Theorem~4.2.5 of~\cite{Prevot2007}),
$
AX_{t_{n}}( \xi) \longrightarrow_{n\rightarrow\infty}AX_{t}( \xi)
$
in $L^{2}( \mathbb{P},\mathfrak{h})$.
Hence Lemma~\ref{lema17} yields
%
%
\begin{equation}
\label{138}
\lim_{n\rightarrow\infty} \mathbb{E}\langle GX_{t_{n}}( \xi)
,AX_{t_{n}}( \xi) \rangle
=
\mathbb{E}\langle GX_{t}( \xi) ,AX_{t}( \xi) \rangle;
\end{equation}
see, for example, Section III.1.7 of~\cite{Kato1995}.
By (\ref{138}) with $A$ replaced by $A ^{*}$,
$t \mapsto\mathbb{E}\langle A ^{*}X_{t}( \xi) ,G X_{t}( \xi) \rangle$
is continuous,
then so is
$t \mapsto\mathbb{E}\langle X_{t}( \xi) ,AG X_{t}( \xi) \rangle$.

We now focus on
$
\sum_{k=1}^{\infty}
\mathbb{E}\langle L_{k}X_{t}( \xi) ,AL_{k}X_{t}( \xi) \rangle
$.
Taking $A=I$ in (\ref{138}) we get
$
\mathbb{E} \Re\langle X_{t_{n}}( \xi) ,GX_{t_{n}}( \xi) \rangle
\rightarrow_{n\rightarrow\infty}
\mathbb{E} \Re\langle X_{t}( \xi) ,GX_{t}( \xi) \rangle$.
Thus
condition (H2.2) leads to
%
%
\begin{equation}
\label{313}
\sum_{k=1}^{\infty}
\mathbb{E}\Vert L_{k}X_{t_{n}}( \xi) \Vert^{2}
\longrightarrow_{n\rightarrow\infty}
\sum_{k=1}^{\infty} \mathbb{E}\Vert L_{k}X_{t}( \xi) \Vert^{2}.
\end{equation}

Using (\ref{313}) we will deduce that
$L_{k}X_{t_{n}}( \xi) $ converges strongly in $L^{2}( \mathbb
{P},\mathfrak{h})$
to $L_{k}X_{t}( \xi) $ as $n \rightarrow\infty$.
Conversely,
suppose that for a given $j \in\mathbb{N}$,
%
%
\begin{equation}
\label{317}
\limsup_{n \rightarrow\infty} \mathbb{E}\Vert L_{j}X_{t_{n}}( \xi)
\Vert^{2}
>
\mathbb{E}\Vert L_{j}X_{t}( \xi) \Vert^{2} .
\end{equation}
Since
$
\mathbb{E}\Vert L_{k}X_{t}( \xi) \Vert^{2}
\leq
\lim\inf_{n\rightarrow\infty}
\mathbb{E}\Vert L_{k}X_{t_{n}}( \xi) \Vert^{2}
$, Fatou's lemma shows
%
%
\begin{equation}
\label{318}
\sum_{k \neq j} \mathbb{E}\Vert L_{k}X_{t}( \xi) \Vert^{2}
\leq
\liminf_{n \rightarrow\infty} \sum_{k \neq j} \mathbb{E}\Vert
L_{k}X_{t_{n}}( \xi) \Vert^{2}.
\end{equation}
According to (\ref{313}) and (\ref{317}) we have
\begin{eqnarray*}
\liminf_{n \rightarrow\infty} \sum_{k \neq j} \mathbb{E}\Vert
L_{k}X_{t_{n}}( \xi) \Vert^{2}
& = &
\sum_{k=1}^{\infty} \mathbb{E}\Vert L_{k}X_{t}( \xi) \Vert^{2}
-
\limsup_{n \rightarrow\infty} \mathbb{E}\Vert L_{j}X_{t_{n}}( \xi)
\Vert^{2} \\
& < &
\sum_{k \neq j} \mathbb{E}\Vert L_{k}X_{t}( \xi) \Vert^{2} ,
\end{eqnarray*}
contrary to (\ref{318}),
and so
%
%
\begin{equation}
\label{332}
\limsup_{n \rightarrow\infty} \mathbb{E}\Vert L_{j}X_{t_{n}}( \xi)
\Vert^{2}
\leq
\mathbb{E}\Vert L_{j}X_{t}( \xi) \Vert^{2}
.
\end{equation}
Applying Lemma~\ref{lema17} we get that
$L_{j}X_{t_{n}}( \xi) $ converges weakly
in $L^{2}( \mathbb{P},\mathfrak{h})$
to $L_{j}X_{t}( \xi) $ as $n \rightarrow\infty$,
and so (\ref{332}) leads to\vspace*{2pt}
$
L_{k}X_{t_{n}}( \xi)
\longrightarrow_{n \rightarrow\infty}
L_{k}X_{t}( \xi)
$
in
$L^{2}( \mathbb{P},\mathfrak{h})$.

From condition (H2.2) it follows that
$
\sum_{k=1}^{n}
\mathbb{E}\langle L_{k}X_{t}( \xi) ,AL_{k}X_{t}( \xi) \rangle
$
converges to
$
\sum_{k=1}^{\infty}
\mathbb{E}\langle L_{k}X_{t}( \xi) ,AL_{k}X_{t}( \xi) \rangle
$
as $n\rightarrow\infty$ uniformly on any finite interval.
Since
$
\mathbb{E}\langle L_{k}X_{t_{n}}( \xi) ,AL_{k}X_{t_{n}}( \xi) \rangle
\longrightarrow_{n\rightarrow\infty}
\mathbb{E}\langle L_{k}X_{t}( \xi) ,AL_{k}X_{t}( \xi) \rangle
$,
the map
$
t \mapsto
\sum_{k=1}^{\infty}
\mathbb{E}\langle L_{k}X_{t}( \xi) ,AL_{k}X_{t}( \xi) \rangle
$
is continuous.
\end{pf}

Third,
we deal with basic properties of
the probabilistic representation of the right-hand side of (\ref{311}).

\begin{lemma}
\label{lema61}
Let Hypothesis~\ref{HipN5} hold.
For any $\xi\in L_{C}^{2}( \mathbb{P},\mathfrak{h}) $,
we define
\[
\mathcal{L}_{*} ( \xi,t )
=
\mathbb{E} \vert G X_{t}( \xi) \rangle\langle X_{t}( \xi) \vert
+
\mathbb{E} \vert X_{t}( \xi) \rangle\langle GX_{t}( \xi) \vert
+
\sum_{k=1}^{\infty}\mathbb{E} \vert L_{k} X_{t}( \xi) \rangle
\langle L_{k}X_{t}( \xi) \vert
.
\]
Then $\mathcal{L}_{*} ( \xi,t )$ is a trace-class operator on
$\mathfrak{h}$ whose trace-norm is uniformly bounded with respect to
$t$ on bounded time intervals;
the series involved in the definition of $\mathcal{L}_{*}$ converges in
$\mathfrak{L} _{1}( \mathfrak{h} )$.
\end{lemma}

\begin{pf}
By condition (H2.2),
using (\ref{55}) and Lemma~\ref{lema51} we get
\begin{eqnarray*}
&& \Vert\mathbb{E} \vert G X_{t}( \xi) \rangle\langle X_{t}( \xi)
\vert\Vert_{1}
+
\Vert\mathbb{E} \vert X_{t}( \xi) \rangle\langle GX_{t}( \xi)
\vert\Vert_{1}
+
\sum_{k=1}^{\infty}
\Vert\mathbb{E} \vert L_{k} X_{t}( \xi) \rangle
\langle L_{k}X_{t}( \xi) \vert
\Vert_{1}
\\
&&\qquad \leq
4 \mathbb{E} ( \Vert X_{t}( \xi) \Vert\Vert GX_{t}( \xi) \Vert)
\leq
K \sqrt{ \mathbb{E} \Vert\xi\Vert^{2} }
\sqrt{ \mathbb{E} \Vert X_{t}( \xi) \Vert_{C}^{2} },
\end{eqnarray*}
where the last inequality follows from
$G \in\mathfrak{L}( ( \mathcal{D}( C) ,\Vert\cdot\Vert_{C})
,\mathfrak{h}) $.
\end{pf}

Applying Lemmata~\ref{lema51} and~\ref{lema18} we easily obtain Lemma
\ref{lema61b}.

\begin{lemma}
\label{lema61b}
Adopt the assumptions and notation of Lemma~\ref{lema61},
together with $A \in\mathfrak{L} ( \mathfrak{h} )$.
Then,
the trace of
$ A \mathcal{L}_{*} ( \xi,t ) $ is equal to
\[
\mathbb{E}\langle X_{t}( \xi) ,AG X_{t}( \xi) \rangle
+
\mathbb{E}\langle GX_{t}( \xi) ,AX_{t}( \xi) \rangle
+
\sum_{k=1}^{\infty}
\mathbb{E}\langle L_{k}X_{t}( \xi) ,AL_{k}X_{t}( \xi) \rangle,
\]
and
$t \mapsto\operatorname{tr} ( A \mathcal{L}_{*} ( \xi,t ) )$
is continuous as a function from $[0, \infty[$
to $ \mathbb{C}$.

\end{lemma}

We proceed to prove that
$\mathbb{E}\vert X_{t}( \xi) \rangle\langle X_{t}( \xi) \vert$
satisfies an integral version of (\ref{3}).
To this end,
we combine the regularity of $X( \xi)$
with It\^{o}'s formula.

\begin{lemma}
\label{lema23}
Adopt Hypothesis~\ref{HipN5} together with $\xi\in L_{C}^{2}( \mathbb
{P},\mathfrak{h}) $.
Then
%
%
\begin{equation}
\label{136}
\rho_{t}( \mathbb{E}\vert\xi\rangle\langle\xi\vert)
=
\mathbb{E}\vert\xi\rangle\langle\xi\vert
+
\int_{0}^{t} \mathcal{L}_{*} ( \xi,s ) \,ds,
\end{equation}
where
$t\geq0$
and
$\mathcal{L}_{*} ( \xi,s )$ is as in Lemma~\ref{lema61};
we understand the above integral in the sense of Bochner integral in
$\mathfrak{L}_{1}( \mathfrak{h}) $.
\end{lemma}

\begin{pf}
Our proof is based on arguments given in Section~\ref{subsecexistence}.
Fix $x \in\mathfrak{h}$,
and
choose
$
\tau_{n}=\inf\{ s\geq0\dvtx\Vert X_{s}( \xi) \Vert>n\}
$,
with $n \in\mathbb{N}$.
Applying the complex It\^{o} formula we obtain that
%
%
\begin{equation}
\label{135}
\langle X_{t\wedge\tau_{n}}( \xi) ,x\rangle X_{t\wedge\tau_{n}}( \xi)
=
\langle\xi,x\rangle\xi+\mathbb{E}\int_{0}^{t\wedge\tau_{n}} L_{x}
(X_{s}( \xi) ) \,ds
+
M_t,
\end{equation}
where
$
M_t
=
\sum_{k=1}^{\infty}\int_{0}^{t\wedge\tau_{n}}( \langle
X_{s}( \xi) ,x\rangle L_{k}X_{s}( \xi)
+\langle L_{k}X_{s}( \xi) ,x\rangle X_{s}(
\xi) ) \,dW_{s}^{k}
$,
and
$
L_{x} (z )
=
\langle z ,x\rangle G z
+ \langle Gz ,x\rangle z
+\sum_{k=1}^{\infty}\langle L_{k} z ,x\rangle L_{k}z
$
for any $z \in\mathcal{D}( C )$.
According to condition (H2.2) we have
\begin{eqnarray*}
&&
\mathbb{E}\sum_{k=1}^{\infty}\int_{0}^{t\wedge\tau_{n}}\Vert\langle
X_{s}( \xi) ,x\rangle L_{k}X_{s}(
\xi) +\langle L_{k}X_{s}( \xi) ,x\rangle X_{s}( \xi) \Vert^{2}\,ds
\\
&&\qquad \leq
4 n^{3} \Vert x\Vert^{2}
\mathbb{E} \int_{0}^{t\wedge\tau_{n}} \Vert G X_{s} \Vert \,ds.
\end{eqnarray*}
Therefore
$
\mathbb{E} M_t =0
$
by $G$ belongs to $\mathfrak{L}( ( \mathcal{D}( C) ,\Vert\cdot\Vert
_{C}) ,\mathfrak{h}) $,
and
so (\ref{135}) yields
%
%
\begin{equation}
\label{131}
\mathbb{E} \langle X_{t\wedge\tau_{n}}( \xi) ,x\rangle X_{t\wedge\tau
_{n}}( \xi)
=
\mathbb{E}\langle\xi,x\rangle\xi+\mathbb{E}\int_{0}^{t\wedge\tau
_{n}} L_{x} (X_{s}( \xi) ) \,ds.
\end{equation}

We will take the limit as $n \rightarrow\infty$ in (\ref{131}).
Since $ X ( \xi)$ has continuous sample paths,
$\tau_{n}\nearrow_{n\rightarrow\infty}\infty$.
By (H2.1)\vspace*{1pt} and (H2.2),
applying the dominated convergence yields
$
\lim_{n \rightarrow\infty} \mathbb{E}\int_{0}^{t\wedge\tau_{n}} L_{x}
(X_{s}( \xi) ) \,ds
=
\mathbb{E}\int_{0}^{t} L_{x} (X_{s}( \xi) ) \,ds.
$
Combinig
$
\mathbb{E} ( \sup_{s \in[0, t+1 ]} \Vert X_s( \xi) \Vert^{2} )
< \infty
$
with the dominated convergence theorem gives
$
\lim_{n \rightarrow\infty} \mathbb{E} \langle X_{t\wedge\tau_{n}}(
\xi) , x\rangle X_{t\wedge\tau_{n}}( \xi)
=
\mathbb{E} \langle X_{t}( \xi) ,x\rangle X_{t}( \xi).
$
Then,
letting first $n \rightarrow\infty$ in (\ref{131}) and then using
Fubini's theorem, we get
%
%
\begin{equation}
\label{63}
\mathbb{E}\langle X_{t}( \xi) ,x\rangle X_{t}( \xi)
=
\mathbb{E}\langle\xi,x\rangle\xi
+
\int_{0}^{t} \mathbb{E} L_{x} (X_{s}( \xi) ).
\end{equation}

By condition (H2.2),
the dominated convergence theorem leads to
\[
\mathbb{E} \sum_{k=1}^{\infty} \langle L_{k} X_{s}( \xi) ,x\rangle
L_{k}X_{s}( \xi)
=
\sum_{k=1}^{\infty}\mathbb{E} \langle L_{k} X_{s}( \xi) ,x\rangle
L_{k}X_{s}( \xi) ,
\]
and so Lemma~\ref{lema51} yields
$
\mathbb{E} L_{x} (X_{s}( \xi) )
=
\mathcal{L}_{*} ( \xi,s ) x
$,
hence
%
%
\begin{equation}
\label{61}
\int_{0}^{t} \mathbb{E} L_{x} (X_{s}( \xi) )
=
\int_{0}^{t} \mathcal{L}_{*} ( \xi,s ) x \,ds .
\end{equation}

Since the dual of $ \mathfrak{L}_{1}( \mathfrak{h})$ consists in all
linear maps
$\varrho\mapsto\operatorname{tr} ( A \varrho)$ with $A \in\mathfrak
{L}( \mathfrak{h})$,
Lemma~\ref{lema61b} implies that
$
t \mapsto\mathcal{L}_{*} ( \xi,t )
$
is measurable as a function from $[0, \infty[$
to $ \mathfrak{L}_{1}( \mathfrak{h})$.
Furthermore,
using Lemma~\ref{lema61} we get that
$
t \mapsto\mathcal{L}_{*} ( \xi,t )
$
is a Bochner integrable $\mathfrak{L}_{1}( \mathfrak{h}) $-valued
function on bounded intervals.
Then
(\ref{63}), together with (\ref{61}),
gives (\ref{136}).
\end{pf}

We are in position to show (\ref{311}) and (\ref{312}) with the help of
Hypothesis~\ref{HipN1}.

\begin{pf*}{Proof of Theorem~\ref{teor10}}
By Theorem~\ref{teorema4},
$
\varrho= \mathbb{E}\vert\xi\rangle\langle\xi\vert
$
for certain $\xi\in L_{C}^{2}( \mathbb{P},\mathfrak{h}) $.
Theorem~\ref{teorema8} now gives
$
A G\rho_{t}(\varrho) = \mathbb{E} \vert A GX_{t} ( \xi) \rangle
\langle X_{t}( \xi) \vert
$.
Applying Hypothesis~\ref{HipN1} we get that
$G^{\ast}, L_{1}^{\ast}, L_{2}^{\ast}, \ldots$ are densely defined
and $G^{\ast\ast}$, $L_{1}^{\ast\ast}, \ldots$ coincide with the
closures of $G, L_{1}, \ldots,$ respectively; see, for example, Theorem
III.5.29 of~\cite{Kato1995}.
Theorem~\ref{teorema8} yields
$
A \rho_{t}(\varrho) G^{\ast} = \mathbb{E} \vert AX_{t}( \xi) \rangle
\langle GX_{t}( \xi) \vert
$
and
$
A L_{k} \rho_{t}(\varrho) L_{k}^{\ast}
=
\mathbb{E} \vert A L_{k}X_{t}( \xi) \rangle\langle L_{k}X_{t}( \xi)
\vert
$.
Therefore
%
%
\begin{equation}
\label{64}
\mathcal{L}_{*} ( \xi,t )
=
G\rho_{t}( \varrho) + \rho_{t}( \varrho)G^{\ast}
+\sum_{k=1}^{\infty} L_{k}\rho_{t}( \varrho) L_{k}^{\ast} ,
\end{equation}
where $\mathcal{L}_{*} ( \xi,t )$ is as in Lemma~\ref{lema61}.
Combining (\ref{64}) with Lemma~\ref{lema23} we get (\ref{311}),
and so
$
\operatorname{tr}( A\rho_{t}( \varrho) )
=
\operatorname{tr}( A \varrho)
+
\int_{0}^{t}
\operatorname{tr}(
A \mathcal{L}_{*} ( \xi,s )
) \,ds
$
for all $t \geq0$.
Using the continuity of $\mathcal{L}_{*} ( \xi,\cdot)$ we obtain (\ref{312}).
\end{pf*}

\subsection{\texorpdfstring{Proof of Theorem \protect\ref{teorema9}}{Proof of Theorem 4.5}}
\label{subsecteorema9}

We first obtain the existence of a solution of (\ref{3}) in the
semigroup sense,
without Hypothesis~\ref{HipN1}.\vadjust{\goodbreak}

\begin{lemma}
\label{lema26}
Under Hypothesis~\ref{HipN5},
$ \rho$
is a semigroup $C$-solution of (\ref{3}).
\end{lemma}

\begin{pf}
By Theorem~\ref{teor8},
$( \rho_{t}) _{t\geq0}$ is a semigroup of bounded operators on
$\mathfrak{L}_{1}( \mathfrak{h}) $
that satisfies property (i) of Definition~\ref{defSemigroupSol}.
Fix $\varrho= \vert x\rangle\langle x\vert$,
with $x\in\mathcal{D}( C) $.
Thus $\varrho$ is a $C$-regular operator,
and so (\ref{310}) leads to property (ii).
Finally,
using Lemmata~\ref{lema61b} and~\ref{lema23} we get
property (iii).
\end{pf}

We next make it legitimate to use in our context the duality relation
between quantum master equations and adjoint quantum master equations.

\begin{lemma}
\label{lema20}
Let Hypothesis~\ref{HipN5} hold.
Suppose that
$A \in\mathfrak{L}( \mathfrak{h}) $
and that
$( \widehat{\rho}_{t}) _{t\geq0}$
is a semigroup $C$-solution of (\ref{3}).
Then $( \widehat{\rho}_{t}^{\ast}( A ) ) _{t\geq0}$ is a $C$-solution
of (\ref{41})
with initial datum $A$,
where $( \widehat{\rho}_{t}^{\ast} ) _{t\geq0}$ is the adjoint
semigroup of $( \widehat{\rho}_{t}) _{t\geq0}$ (see, e.g.,~\cite{Pazy1983}),
that is,
$( \widehat{\rho}_{t}^{\ast})
_{t\geq0} $ is the unique semigroup of bounded operators on $\mathfrak
{L}( \mathfrak{h}) $ such that
for all $B\in\mathfrak{L}( \mathfrak{h}) $ and $\varrho\in\mathfrak
{L}_{1}( \mathfrak{h}) $,
%
%
\begin{equation}
\label{315}
\operatorname{tr}( \widehat{\rho}_{t}( \varrho) B)
=
\operatorname{tr}( \widehat{\rho}_{t}^{\ast}( B) \varrho) .
\end{equation}
\end{lemma}

\begin{pf}
Using (\ref{315}) we get that for all vectors $x,y \in\mathfrak{h}$
whose norm is~$1$,
\begin{eqnarray*}
\vert\langle y,\widehat{\rho}_{t}^{\ast}( A) x\rangle\vert
& = &
\vert{\operatorname{tr}}( \widehat{\rho}_{t}^{\ast}( A) \vert x\rangle
\langle y\vert) \vert =
\operatorname{tr}( \vert\widehat{\rho}_{t}( \vert x\rangle\langle
y\vert) A\vert)\\
&\leq&
\Vert A\Vert\Vert\widehat{\rho}_{t}\Vert_{\mathfrak{L}( \mathfrak
{L}_{1}( \mathfrak{h} ) )}
\operatorname{tr}( \vert\vert x\rangle\langle y\vert\vert).
\end{eqnarray*}
We conclude from (\ref{55}) that
$ \operatorname{tr}( \vert\vert x\rangle\langle y\vert\vert) =1$,
hence that
$
\vert\langle y,\widehat{\rho}_{t}^{\ast}( A) x\rangle\vert
\leq
\Vert A\Vert\times\Vert\widehat{\rho}_{t}\Vert_{\mathfrak{L}( \mathfrak
{L}_{1}( \mathfrak{h}) ) }
$,
and finally that
%
%
\begin{equation}\label{316}
\Vert\widehat{\rho}_{t}^{\ast}( A) \Vert_{\mathfrak{L}( \mathfrak{h}
) }
\leq\Vert A\Vert\Vert\widehat{\rho}_{t}\Vert_{\mathfrak{L}( \mathfrak
{L}_{1}( \mathfrak{h}) ) }.
\end{equation}
Applying property (i) of Definition~\ref{defSemigroupSol} gives
property (b) of Definition~\ref{definicion3}.

In order to verify property (a),
we will prove the continuity of
$t\mapsto\langle x,\break\widehat{\rho}_{t}^{\ast}( A) y \rangle$
for any $x ,y \in\mathfrak{h}$.
As in the proof of Lemma~\ref{lema42},
we define
$R_{n} = n (n+C )^{-1}$ for $n \in\mathbb{N}$.
According to (\ref{315}) we have
\[
\langle R_{n}x,\widehat{\rho}_{t}^{\ast}( A) R_{n} x\rangle
=
\operatorname{tr}( \widehat{\rho}_{t}^{\ast} ( A ) \vert R_{n} x \rangle
\langle R_{n} x\vert)
=
\operatorname{tr}( \widehat{\rho}_{t} ( \vert R_{n} x \rangle\langle
R_{n} x\vert)
A) .
\]
Since $ R_{n} x\in\mathcal{D}( C) $,
property (ii) of Definition~\ref{defSemigroupSol} implies the
continuity of the function
$t\mapsto\langle R_{n} x,\widehat{\rho}_{t}^{\ast}( A) R_{n} x\rangle$.
By (\ref{316}),
\begin{eqnarray*}
&&
\vert\langle x,\widehat{\rho}_{t}^{\ast}( A) x\rangle-\langle
x,\widehat{\rho}_{s}^{\ast}( A)
x\rangle\vert\\
&&\qquad\leq
\vert\langle R_{n} x,\widehat{\rho}_{t}^{\ast}( A)
R_{n} x\rangle-\langle R_{n} x,\widehat{\rho}_{s}^{\ast}(
A) R_{n} x\rangle\vert
\\
&&\qquad\quad{}+ 2 \Vert A\Vert\bigl( \Vert\widehat{\rho}_{t}\Vert_{\mathfrak{L}(
\mathfrak{L}_{1}( \mathfrak{h}) ) }
+
\Vert\widehat{\rho}_{s}\Vert_{\mathfrak{L} ( \mathfrak{L}_{1}(
\mathfrak{h}) ) }\bigr)
\Vert x\Vert\Vert x - R_{n} x\Vert.
\end{eqnarray*}
Using $R_{n} x \longrightarrow_{n \rightarrow\infty} x$ we deduce that
the map
$t\mapsto\langle x,\widehat{\rho}_{t}^{\ast}( A) x\rangle$ is continuous,
so is
$t\mapsto\langle x,\widehat{\rho}_{t}^{\ast}( A) y \rangle$
by the polarization identity.

Assume that $x \in\mathcal{D}( C )$.
By (\ref{315}), combining
$
\widehat{\rho}_{t+s}^{\ast}( A)
=
\widehat{\rho}_{s}^{\ast}( \widehat{\rho}_{t}^{\ast}( A))
$
with property (iii) of Definition~\ref{defSemigroupSol} yields
\begin{eqnarray*}
&& \lim_{s\rightarrow0+}\frac{1}{s}\bigl( \langle x,\widehat{\rho
}_{t+s}^{\ast}( A) x\rangle-\langle x,\widehat{\rho}_{t}^{\ast}( A)
x\rangle\bigr) \\[-2pt]
&&\qquad =
\lim_{s\rightarrow0+}\frac{1}{s}\bigl( \operatorname{tr} ( \widehat{\rho
}_{s}( \vert x\rangle\langle x\vert) \widehat{\rho}_{t}^{\ast}( A) )
-\operatorname{tr}( \vert x\rangle\langle x\vert\widehat{\rho
}_{t}^{\ast}( A) ) \bigr)\\[-2pt]
&&\qquad=
\mathcal{L} ( \widehat{\rho}_{t}^{\ast}( A),x )
\end{eqnarray*}
with
$
\mathcal{L} ( \widehat{\rho}_{t}^{\ast}( A),x )
=
\langle x,\widehat{\rho}_{t}^{\ast}( A)
Gx\rangle+\langle Gx,\widehat{\rho}_{t}^{\ast}( A)
x\rangle+\sum_{k=1}^{\infty}\langle L_{k}x,\widehat{\rho}%
_{t}^{\ast}( A) L_{k}x\rangle
$.
Thus
%
%
\begin{equation}
\label{325}
\frac{d}{dt}^{+}\langle x,\widehat{\rho}_{t}^{\ast}( A)
x\rangle
=
\mathcal{L} ( \widehat{\rho}_{t}^{\ast}( A),x ).
\end{equation}
From\vspace*{1pt} (\ref{316}) and condition (H2.2)
we get that
$\sum_{k=1}^{\infty}\langle
L_{k}x,\widehat{\rho}_{t}^{\ast}( A) L_{k}x\rangle$ is uniformly
convergent on bounded intervals,
and so
$
t\mapsto\sum_{k=1}^{\infty}\langle L_{k}x,\widehat{\rho}_{t}^{\ast}(
A) L_{k}x\rangle
$
is continuous, and
hence the application
$t\mapsto\frac{d}{dt}^{+}\langle x,\widehat{\rho}_{t}^{\ast}( A)
x\rangle$ is continuous.\vspace*{1pt}
Therefore $\langle x,\widehat{\rho}_{t}^{\ast}( A) x\rangle$
is continuously differentiable (see, e.g., Section 2.1 of~\cite{Pazy1983}).
Property (a) of Definition~\ref{definicion3} now follows from (\ref{325}).
\end{pf}

We are in position to show our second main theorem.

\begin{pf*}{Proof of Theorem~\ref{teorema9}}
Let $( \widehat{\rho}_{t})_{t\geq0}$ be a semigroup $C$-solution of
(\ref{3}). Consider the adjoint semigroup $ (
\widehat{\rho}_{t}^{\ast} ) _{t\geq0}$ of $( \widehat{\rho}_{t})
_{t\geq0}$, and let $( \mathcal{T}_{t}( A) ) _{t\geq0}$ be given by
Theorem~\ref{teorema3}. Combining Lemma~\ref{lema20} with Theorem
\ref{teorema3} we obtain $\widehat{\rho}_{t}^{\ast}( A) =
\mathcal{T}_{t}( A) $ for all $t\geq0$ and $A\in\mathfrak{L}(
\mathfrak{h}) $. If $\varrho\in\mathfrak{L}_{1,C}^{+}( \mathfrak{h}) $
and $A\in\mathfrak{L}( \mathfrak{h}) $, then applying (\ref{315}) and
Lemma~\ref{lema10} yields
\[
\operatorname{tr}( \rho_{t}( \varrho) A)
= \operatorname{tr}( \mathcal{T}_{t}( A) \varrho)
= \operatorname{tr}( \widehat{\rho}_{t}^{\ast}( A) \varrho)
= \operatorname{tr}( \widehat{\rho}_{t}( \varrho) A)
\]
and so
$\rho_{t}( \varrho)
= \widehat{\rho}_{t}( \varrho) $.
Lemma~\ref{lema12} now implies that
$\rho_{t}( \varrho)
= \widehat{\rho}_{t}( \varrho) $
for all $\varrho$ belonging to $\mathfrak{L}_{1}^{+}( \mathfrak{h}) $, hence
$\rho_{t} = \widehat{\rho}_{t}$.
Finally, Lemma~\ref{lema26} completes the proof.
\end{pf*}

\subsection{\texorpdfstring{Proof of Theorem \protect\ref{teorema5}}{Proof of Theorem 4.6}}
\label{subsecteorema5}

From~\cite{FagnolaMora2010}
we have that Hypothesis~\ref{HipN5} holds with $C = P^2 + Q^2$.
Hence Theorem~\ref{teorema9} yields our first assertion.

Suppose that $A=P$ or $A=Q$.
Using, for instance, the spectral theorem,
we deduce the existence of a sequence $A_n$
of bounded self-adjoint operators in $\mathfrak{h}$
such that for all $f \in\mathcal{D} ( A )$ we have
$\| A_n f \| \leq\| A f \|$
and
$A_n f \longrightarrow_{n \rightarrow\infty} Af$.
Applying Theorems~\ref{teorema8} and~\ref{teor10}
(or better Lemmata~\ref{lema61b} and~\ref{lema23})
gives
%
%
\begin{eqnarray}
\label{32}\quad
&&\operatorname{tr}( A_n \rho_{t}( \varrho) )\nonumber\\
&&\qquad=
\operatorname{tr}( A_n \varrho)
+
\int_{0}^{t}
\Biggl(
\mathbb{E}\langle A_n X_{t}( \xi) , G X_{t}( \xi) \rangle
+
\mathbb{E}\langle GX_{t}( \xi) ,A_nX_{t}( \xi) \rangle\\
&&\qquad\quad\hspace*{65.3pt}\hspace*{52.4pt}{}
+
\sum_{k=1}^{\infty}
\mathbb{E}\langle L_{k}X_{t}( \xi) ,A_nL_{k}X_{t}( \xi) \rangle
\,ds
\Biggr),\nonumber
\end{eqnarray}
where
$
\varrho= \mathbb{E} \vert\xi\rangle\langle\xi\vert
$
with
$
\mathbb{E} ( \Vert C \xi\Vert^{2} + \Vert\xi\Vert^{2} )
<\infty
$.
By the dominated convergence theorem,
letting $n \rightarrow\infty$ we obtain
%
%
\begin{eqnarray}
\label{32n}\quad
&&\operatorname{tr}( A \rho_{t}( \varrho) )\nonumber\\
&&\qquad=
\operatorname{tr}( A \varrho)
+
\int_{0}^{t}
\Biggl(
\mathbb{E}\langle A X_{t}( \xi) , G X_{t}( \xi) \rangle
+
\mathbb{E}\langle GX_{t}( \xi) ,A X_{t}( \xi) \rangle
\\
&&\hspace*{140pt}{}+
\sum_{k=1}^{\infty}
\mathbb{E}\langle L_{k}X_{t}( \xi) ,A L_{k}X_{t}( \xi) \rangle
\,ds
\Biggr).\nonumber
\end{eqnarray}

Let $f \in C_{c}^{\infty} ( \mathbb{R}, \mathbb{C} )$.
Using
$
[ P, Q ] = -i I
$
we get that
$ \mathcal{L} ( P ) f = i [ H, P ] f $
and
$ \mathcal{L} ( Q ) f = i [ H, Q ] f $.
Therefore
%
%
\begin{eqnarray}
\label{34}
&&\langle A f, Gf \rangle
+
\langle Gf, Af \rangle
+\sum_{k=1}^{\infty}\langle L_{k}f, A L_{k}f
\rangle\nonumber\\[-8pt]\\[-8pt]
&&\qquad=
\cases{
\langle f, P f \rangle/ m ,
&\quad
if $A=Q$,\cr
-2c \langle f, Q f \rangle, &\quad
if $A=P$.}\nonumber
\end{eqnarray}
Since $C_{c}^{\infty} ( \mathbb{R}, \mathbb{C} )$ is a core for
$C = P^2 + Q^2$,
combining a limit procedure with, for instance, Lemma 12 of \cite
{FagnolaMora2010}
we get that (\ref{34}) holds for all $f \in\mathcal{D} ( C )$.
Then,
(\ref{32n}) leads to (\ref{31n}).

\subsection{\texorpdfstring{Proof of Theorem \protect\ref{corolario2}}{Proof of Theorem 5.1}}
\label{subseccorolario2}

Let $\xi$ be distributed according to $\theta$.
Set $\widetilde{\mathbb{Q}} = \Vert X_{T}( \xi) \Vert^{2}\cdot\mathbb{P}$,
with $T > 0$.
For any $t \in[ 0, T ]$,
we choose
$ \widetilde{Y}_{t} = X_{t}( \xi) /\Vert X_{t}( \xi) \Vert$
if $X_{t}( \xi) \neq0$
and
$ \widetilde{Y}_{t} = 0 $ otherwise;
let
$
B_{t}^{k}=W_{t}^{k}-\int_{0}^{t}\frac{1}{\Vert X_{s}(
\xi) \Vert^{2}}\,d[ W^{k}, X( \xi) ]_{s}
$
for any $k\in\mathbb{N}$.
Proceeding along the same lines as in the proof of Proposition~1 of
\cite{MoraReAAP2008}
we obtain\vspace*{1pt} that
$
( \mathbb{Q},( Y_{t}) _{t\in[ 0,T] },( B_{t}^{k}) _{t\in[ 0,T] }^{k\in
\mathbb{N}})
$
is a $C$-solution of (\ref{5}) with initial law~$\theta$.
By Remark~\ref{nota6},
(\ref{5}) has a unique $C$-solution with initial distribution $\theta$.
Therefore
the distribution of $\widetilde{Y}_{t}$ with respect to
$\widetilde{\mathbb{Q}}$ coincides with the distribution of $Y_{t}$
under $\mathbb{Q}$.
From~\cite{FagnolaMora2010} we have that
$( \Vert X_{t}\Vert^{2})_{t\in[ 0,T] }$ is a martingale, and
hence
for any $x\in\mathfrak{h}$
and $t \in[ 0, T ]$,
\[
\mathbb{E}_{\mathbb{Q}}\vert\langle x,Y_{t}\rangle\vert^{2}
=
\mathbb{E}_{\widetilde{ \mathbb{Q}}}\vert\langle x,\widetilde
{Y}_{t}\rangle\vert^{2}
=
\mathbb{E}_{\mathbb{P}}( \vert\langle x,\widetilde{Y}_{t}\rangle\vert
^{2}\Vert X_{t}( \xi)
\Vert^{2})
=
\mathbb{E}_{\mathbb{P}}\vert\langle x,X_{t}( \xi) \rangle\vert^{2}.
\]
Applying (\ref{31}) and the polarization identity
gives $\rho_{t}( \varrho) =\mathbb{E} \vert Y_{t}\rangle\langle
Y_{t}\vert$.

\subsection{\texorpdfstring{Proof of Theorem \protect\ref{teorema7}}{Proof of Theorem 5.2}}
\label{subsecteorema7}

%
Let $( \mathbb{Q},( Y_{t}) _{t\geq0},(
B_{t})_{t\geq0}) $ be the $C$-solution of (\ref{5}) with
initial distribution $\Gamma$; see Remark~\ref{nota6}. Choose
$\varrho_{\infty}=\mathbb{E}\vert Y_{0}\rangle\langle
Y_{0}\vert$. Then, Theorem~\ref{corolario2} shows that $
\rho_{t}( \varrho_{\infty})
=
\mathbb{E}\vert Y_{t}\rangle\langle Y_{t}\vert
$
for all $t\geq0$.

As in the proof of Theorem 3 of~\cite{MoraReAAP2008},
applying techniques of well-posed martigale problems
we obtain the Markov property of the $C$-solutions of (\ref{5})
under Hypothesis~\ref{HipN5}.
Hence for any $x\in\mathfrak{h}$ and $t\geq0$
\[
\mathbb{E}\bigl( \mathbf{1}_{[ 0,\Vert x\Vert^{2}]
}( \vert\langle x,Y_{t}\rangle\vert^{2})\bigr)=
\mathbb{E}\biggl( \int_{\mathfrak{h}}\mathbf{1}_{[ 0,\Vert
x\Vert^{2}] }( \vert\langle x,y\rangle
\vert^{2}) P_{t}( Y_{0},dy) \biggr).\vadjust{\goodbreak}
\]
On the other hand,
using (\ref{IM1}) we deduce that
\[
\mathbb{E}\bigl( \mathbf{1}_{[ 0,\Vert x\Vert^{2}]
}( \vert\langle x,Y_{0}\rangle\vert^{2})
\bigr)
=
\int_{\mathfrak{h}}\biggl( \int_{\mathfrak{h}}\mathbf{1}_{[
0,\Vert x\Vert^{2}] }( \vert\langle
x,y\rangle\vert^{2}) P_{t}( z,dy) \biggr)
\Gamma( dz) .
\]
Now, combining $ \Vert Y_{t} \Vert= 1 $ with $ \int_{\mathfrak{h}}(
\int_{\mathfrak{h}}\mathbf{1}_{[ 0,\Vert x\Vert^{2}] }( \vert\langle
x,y\rangle\vert^{2}) P_{t}( z,dy) ) \Gamma( dz) =\break \mathbb{E}(
\int_{\mathfrak{h}}\mathbf{1}_{[ 0,\Vert x\Vert^{2}] }( \vert\langle
x,y\rangle \vert^{2}) P_{t}( Y_{0},dy) ) $ we get\vspace*{1pt} $
\mathbb{E}\vert\langle x,Y_{0}\rangle\vert ^{2}=\mathbb{E}\vert\langle
x,Y_{t}\rangle\vert ^{2} $. This gives $ \mathbb{E}\vert
Y_{t}\rangle\langle Y_{t}\vert = \mathbb{E}\vert Y_{0}\rangle\langle
Y_{0}\vert $ and so $\rho_{t}( \varrho _{\infty}) =\varrho_{\infty}$.

\subsection{\texorpdfstring{Proof of Theorem \protect\ref{teor14}}{Proof of Theorem 6.1}}
\label{subsecteor14}

Since $\mathcal{D}( G ) = \mathcal{D}( N^{4} )$, from Remark
\ref{nota3} we have that $G$ is a closable operator satisfying $G
\in\mathfrak{L}( ( \mathcal{D}( N^{p}), \mbox{$\Vert\cdot\Vert_{ N^{p}}$})
,\mathfrak{h})$. Fix $x \in\mathfrak{h}$ such that $x_{n} := \langle
e_{n}, x \rangle$ is equal to $0$ for all $n \in\mathbb{Z}_{+}$ except
a finite number. An easy computation shows that $ 2\Re\langle N^{2p} x,
Gx\rangle+\sum_{k=1}^{\infty} \Vert N^{p} L_{k}x \Vert^{2} $ is equal
to the sum of $ 4 p ( \vert\alpha_{5} \vert^{2} - \vert\alpha_{4}
\vert^{2} ) \sum _{n = 0}^{\infty} n^{2p+1} \vert x_{n} \vert^{2} $ and
\[
2 \beta_{1} \sum_{n=1}^{\infty} \sqrt{n+1} \bigl( ( n+1 )^{2p} - n^{2p} \bigr)
\Re( x_{n} \overline{ x_{n+1}} )
+
\sum_{n = 0}^{\infty} f(n ) \vert x_{n} \vert^{2} ,
\]
where $f$ is a $2p$-degree polynomial whose coefficients depend on $
\vert\alpha_{k} \vert^{2} $ with $k=1,2,4,5$. Hence
$N^{p}$ satisfies Hypothesis~\ref{HipN4} whenever
$\vert\alpha_{4}\vert\geq\vert\alpha_{5}\vert$. From
\cite{FagnolaMora2010} it follows that $N^{p}$ fulfills condition
(H2.3) of Hypothesis~\ref{HipN5}, and so Theorems~\ref{teor10} and
\ref{teorema9} lead to statement (i).

From Theorem 8 of~\cite{MoraReAAP2008} we have the existence of an
invariant probability measure $\Gamma$ for (\ref{5}) that satisfies the
properties given in Hypothesis~\ref{Hip2} with $C=N^p$. Using Theorem
\ref{teorema7} yields statement (ii).

\section*{Acknowledgments}

The author wishes to express his gratitude to the anonymous referees,
whose suggestions and constructive criticisms led to substantial
improvements in the presentation. Moreover, I thank Franco Fagnola and
Roberto Quezada for helpful comments.


%

\printaddresses

\end{document}